\documentclass[12pt]{article}
\usepackage[margin=1in]{geometry}
\usepackage{amsmath,amssymb,amsthm,mathtools,booktabs,array}
\usepackage[authoryear,round]{natbib}
\usepackage{setspace}
\usepackage[colorlinks=true,citecolor=blue,urlcolor=blue,linkcolor=blue]{hyperref}
\theoremstyle{plain}
\newtheorem{theorem}{Theorem}
\newtheorem{lemma}{Lemma}
\newtheorem{proposition}{Proposition}
\newtheorem{corollary}{Corollary}
\numberwithin{equation}{section}
\makeatletter
\renewcommand\section{\@startsection{section}{1}{\z@}{-3.0ex plus -1ex minus -.2ex}{1.5ex plus .2ex}{\normalfont\large\bfseries}}
\renewcommand\subsection{\@startsection{subsection}{2}{\z@}{-2.5ex plus -1ex minus -.2ex}{1.0ex plus .2ex}{\normalfont\normalsize\bfseries}}
\makeatother
\newcommand{\R}{\mathbb R}
\newcommand{\tr}{\operatorname{tr}}
\newcommand{\diag}{\operatorname{diag}}
\newcommand{\rank}{\operatorname{rank}}
\newcommand{\pos}{\operatorname{pos}}
\newcommand{\rev}[1]{\widetilde{#1}}
\newcommand{\Deff}{D_{\mathrm{eff}}}
\newcommand{\lmin}{\lambda_{\min}}
\begin{document}
\begin{center}
{\large\bfseries Exact optimality and foldover designs in order-of-addition experiments}\\[14pt]
Stelios D. Georgiou$^{1}$, Haralambos Evangelaras$^{2}$ and Stella Stylianou$^{1}$\\[8pt]
{\itshape $^{1}$RMIT University and $^{2}$University of Piraeus}
\end{center}

\begin{quotation}
\noindent{\itshape Abstract:} In an order-of-addition experiment the response depends on the order in which several components are added. Under the pair-wise ordering model, order-of-addition orthogonal arrays are D-optimal, but with four or more components they exist only when the run size is a multiple of twelve. This paper gives exact optimality results for other run sizes. We show that adding any one order to an order-of-addition orthogonal array yields a design that is D-optimal among all designs with one more run, and that is also optimal under the I-criterion and two further criteria, so the added order can be chosen freely. For foldover designs, in which every order is accompanied by its reverse, adding any foldover pair to a foldover order-of-addition orthogonal array is optimal among foldover designs, and it is optimal among all designs when the mean contains an arbitrary term that does not change when an order is reversed. The proofs use a lower bound on the distance between the moment matrix of a design and that of the full design, and for foldover designs with four components this distance is minimised exactly at every run size. Deleting is not the counterpart of adding: examples with five, six and seven components show that an orthogonal array with one run deleted, or a foldover array with one pair deleted, need not be D-optimal. We also classify the D-optimal foldover designs with five components and twenty runs, and compare foldover designs with other designs when effects outside the model are present.

\vspace{9pt}
\end{quotation}

\section{Introduction}\label{sec:intro}

In an order-of-addition (OofA) experiment, $m$ components are added one at a time and the response depends on the order in which they are added. The order in which drugs are given in a combination therapy is a well-known example \citep{ding2015,wangxuding2020}, and many others from chemistry, food science and scheduling are described in the reviews of \citet{linpeng2019} and \citet{linrios2025}. There are $m!$ possible orders. The full design, which uses each of them once, is rarely affordable when $m$ exceeds four or five, and a fraction has to be chosen.

The model used most often for these experiments is the pair-wise ordering (PWO) model of \citet{vannostrand1995}. Its regressors are the $q=m(m-1)/2$ PWO factors, each of which records which of two components is added first. Under this model \citet{peng2019} showed that the full design is optimal for a broad class of criteria. \citet{voelkel2019} introduced order-of-addition orthogonal arrays (OofA-OAs), which are fractions that keep the pairwise balance of the PWO factors in the full design. An OofA-OA has the same moment matrix as the full design, a fraction is as efficient as the full design under the $D$-criterion if and only if it is an OofA-OA \citep{zhao2022}, and such a design is also $I$-optimal \citep{schoenmee2023}. Such arrays are constructed or enumerated by \citet{peng2019}, \citet{chen2020}, \citet{tsai2022}, \citet{zhaodong2022} and \citet{schoenmee2023}.

For $m\geq4$ the run size of an OofA-OA must be a multiple of 12 \citep{voelkel2019}. Other run sizes arise all the time: one or two further runs can be afforded, a run is lost, or the resources simply do not come in twelves. \citet{schoenmee2023} suggest adding a few runs to an OofA-OA with optimal design software, or deleting its runs one at a time, and designs for such run sizes have also been found by computer search or by combining smaller designs \citep{voelkel2019,winker2020,zhao2021,tsai2022}. These designs are judged by their $D$-efficiency relative to the full design. The full design is an unattainable benchmark at these run sizes, however, so an efficiency below one does not say how far a design is from the best design of its size. For three and four components the best designs with 4 and 7 runs are known from the complete search of \citet{zhao2021}. The one general statement of optimality at such run sizes that we know of is Theorem~2 of \citet{tsai2022}, according to which an OofA-OA with one run removed is $D$-optimal among all designs with one run fewer. It is true in the cases with three and four components that we could enumerate, but not in general. For five components and 11 runs the designs found by \citet{winker2020} and \citet{wangwang2023} already have a larger determinant, and Section~\ref{sec:deletion} gives designs with six and seven components for which the same holds.

For an added run the picture is different. We show that an OofA-OA with one added order is $D$-optimal among all designs of the enlarged run size, whichever order is added and whether or not the array already contains it (Theorem~\ref{thm:mainaugmentation}). The same design minimises the average variance of prediction and is optimal under two further criteria. An experimenter who can afford 25 runs with five or seven components may therefore take a 24-run OofA-OA and add the order of greatest interest, for instance the order in current use, and no 25-run design is better.

Our second group of results concerns foldover designs, in which every order is run together with its reverse. \citet{mee2020} showed that in a design made of such foldover pairs the estimators of the PWO effects are not biased by two-factor interactions of the PWO factors, and foldover designs are studied further by \citet{wangmee2022}, \citet{tsai2025} and \citet{shao2026}. We show that adding any foldover pair to a foldover OofA-OA gives a design that is optimal among all foldover designs with two more runs. Under the PWO model this restriction to foldover designs has a price, which we compute. The picture changes when the mean contains an arbitrary term that is unchanged by reversing the order: the same design is then optimal among all designs (Theorem~\ref{thm:foldoveraug}). The reason is a general fact, Theorem~\ref{thm:nuisancecompleteclass}, which says that under such a model every $D$-optimal design with an even number of runs consists of foldover pairs. The fact does not depend on the PWO factors, and for two-level factorial experiments it turns the foldover principle of \citet{boxhunter1961} into a statement of optimality.

Results of this type are classical for two-level factorial designs. The optimality of orthogonal arrays with one added run was established, under suitable conditions, by \citet{mukerjee1999} and \citet{chai2002}, additions and deletions of several runs are treated by \citet{butlerramos2007}, and $D$-optimal designs whose run sizes are not multiples of four go back to \citet{ehlich1964}, \citet{cheng1980} and \citet{galilkiefer1980}; recent work includes \citet{wangyue2022} and \citet{hameed2025}. Their arguments do not carry over to the PWO model. The PWO factors cannot be set independently, since their levels in a run must agree with a single order, and the moment matrix of the full design is not the identity matrix, so that conditions on the integer entries of $X^\top X$ are of little direct use. We work with the moment matrix of a design on the scale of the full design and measure its distance from the identity. Building on an identity of \citet{mee2020}, we represent the PWO factors by triples of components (Lemma~\ref{lem:triples}) and write this distance as a sum of squares of numbers that are integers or integers shifted by a third of the run size. Parity and divisibility by three then give a sharp lower bound for odd run sizes. For foldover designs with four components we find the minimum distance, and every design attaining it, at each run size (Theorem~\ref{thm:m4discrepancyall}); the answer depends only on the number of foldover pairs modulo 12. This case in turn gives bounds for every $m$ and every run size.

The same tools explain why deleting a run is not the counterpart of adding one. An OofA-OA with one run deleted is as close to the full design as a design of its size can be, but among the moment matrices at that distance it has the smallest determinant, whereas an added run gives the largest. Foldover pairs behave in the same way, and for seven components and 46 runs we exhibit a foldover design with a larger determinant than a foldover OofA-OA with a pair deleted. The simple picture also ends at run sizes that are not next to a multiple of 12. With five components and 20 runs no OofA-OA exists, and a complete search finds exactly five nonequivalent $D$-optimal foldover designs, which share their determinant but differ under the other criteria (Theorem~\ref{thm:global5}). A foldover design, moreover, protects only against effects that are unchanged by reversal, and a comparison under two families of omitted effects indicates when this protection outweighs the loss of efficiency.

Section~\ref{sec:model} gives the notation, Section~\ref{sec:adding} the result on added runs, and Section~\ref{sec:bound} the bounds and the proof of Theorem~\ref{thm:mainaugmentation}. Section~\ref{sec:deletion} turns to deleted runs, Section~\ref{sec:consequences} treats foldover designs, and Section~\ref{sec:discussion} concludes. Proofs that are not given in the text, the designs and the enumerations are in the Supplementary Material.

\section{Preliminaries}\label{sec:model}

\subsection{The PWO model and orthogonal arrays}

Let $\Omega_m$ be the set of the $m!$ permutations of $\{1,\ldots,m\}$, where $m\geq3$. A permutation $a=(a_1,\ldots,a_m)$ is an order of addition in which component $a_1$ is added first and $a_m$ last. Let $\pos_a(j)$ be the position of component $j$ in $a$. For $j<k$ the PWO factor is
\[
z_{jk}(a)=\begin{cases}
 1,&\pos_a(j)<\pos_a(k),\\
 -1,&\pos_a(j)>\pos_a(k),
\end{cases}
\]
and $z(a)$ is the $q\times1$ vector of the PWO factors, with the pairs in lexicographic order and $q=m(m-1)/2$. The PWO model is
\[
E\{Y(a)\}=\beta_0+z(a)^\top\beta,
\]
with $p=q+1$ parameters and uncorrelated errors of mean zero and variance $\sigma^2$. A design $D$ with $n$ runs is a collection of $n$ orders in which an order may appear more than once. Let $Z_D$ be its $n\times q$ matrix of PWO factors, and write $x(a)=(1,z(a)^\top)^\top$, $X_D=[\mathbf1\ Z_D]$ and
\[
M_D=X_D^\top X_D/n
\]
for the moment matrix of $D$. We shall also need design measures. A design measure puts weights $w_a\geq0$ with $\sum_aw_a=1$ on the orders and has moment matrix $\sum_aw_ax(a)x(a)^\top$, and a design with $n$ runs corresponds to weights proportional to its numbers of replications.

Let $B$ be the $m\times q$ matrix whose column for the pair $j<k$ is $e_j-e_k$, where $e_j$ is the $j$th unit vector of $\R^m$. Since $BB^\top=mI_m-\mathbf1\mathbf1^\top$, the full design $F$, which contains every order of $\Omega_m$ once, has
\begin{equation}\label{eq:full}
 M_F=\diag(1,\Sigma),\qquad
 \Sigma=\frac{I_q+B^\top B}{3},\qquad
 \det\Sigma=\frac{(m+1)^{m-1}}{3^q}
\end{equation}
\citep{peng2019}. Following \citet{voelkel2019}, a design is an OofA-OA of strength $t$ if in every $t$ columns of $Z_D$ the level combinations occur in the same proportions as in the full design. In this paper an OofA-OA has strength 2 unless we say otherwise. Since the PWO factors have two levels, $D$ is an OofA-OA if and only if $M_D=M_F$, and this is the case if and only if $D$ attains the determinant of the full design \citep[Theorem~4]{zhao2022}. The run size of an OofA-OA is a multiple of 12 when $m\geq4$ \citep{voelkel2019}. It is a multiple of 6 when $m=3$, because for three components the moments of order at most two determine the frequencies of the six orders. The $D$-efficiency of a design is
\[
\Deff(D)=\left\{\frac{\det M_D}{\det M_F}\right\}^{1/p},
\]
which is zero for a singular design.

\subsection{Foldover designs}

The reverse of an order $a$ is $\rev a=(a_m,\ldots,a_1)$, and $z(\rev a)=-z(a)$. As in \citet{mee2020} we call $\{a,\rev a\}$ a foldover pair. A foldover design contains $a$ and $\rev a$ equally often, for every order $a$. It can be written $D=H\cup\rev H$, where the half-design $H$ contains one order from each foldover pair of $D$ and a pair that occurs several times in $D$ is counted each time. If $H$ has $h$ orders, then
\begin{equation}\label{eq:fold}
 n=2h,\qquad M_D=\diag(1,G_D),\qquad
 G_D=\frac1h Z_H^\top Z_H,\qquad
 \rank X_D=1+\rank Z_H.
\end{equation}
Replacing an order of $H$ by its reverse changes the sign of a row of $Z_H$ and leaves $G_D$ unchanged. All PWO effects are estimable only if $h\geq q$ and $\rank Z_H=q$. A foldover design with $G_D=\Sigma$ is a foldover OofA-OA. Two foldover designs are equivalent if one is obtained from the other by relabelling the components, the order of the runs being ignored.

\subsection{Criteria}\label{sec:criteria}

For a design $D$ let $\mathcal C_D=M_F^{-1/2}M_DM_F^{-1/2}$ be its moment matrix on the scale of the full design. Besides the determinant we use three criteria. The $I$-criterion is the average variance of prediction over all orders,
\[
 I(D)=\frac1{m!}\sum_{a\in\Omega_m}x(a)^\top M_D^{-1}x(a)=\tr(M_D^{-1}M_F),
\]
so that $\sigma^2I(D)/n$ is the average variance of the fitted mean, and we put $I(D)=+\infty$ for a singular design. The minimum relative efficiency $\lmin(\mathcal C_D)$ compares the precision with which $D$ and the full design, with the same number of runs, estimate a linear function $\ell^\top(\beta_0,\beta^\top)^\top$, for the least favourable $\ell$. The squared distance $\tr\{(\mathcal C_D-I_p)^2\}$ measures how far $M_D$ is from $M_F$. Every design has $\tr\mathcal C_D=p$, as will be seen below, so minimising this distance is the $(M,S)$-criterion of \citet{eccleston1974} applied to $\mathcal C_D$.

For a foldover design $\mathcal C_D=\diag(1,C_D)$ with $C_D=\Sigma^{-1/2}G_D\Sigma^{-1/2}$, and we write
\begin{equation}\label{eq:discrepancy}
 \Delta(D)=\tr\{(C_D-I_q)^2\},\qquad
 \theta(D)=\sqrt{\frac{\Delta(D)}{q(q-1)}}.
\end{equation}
The matrix $C_D$ is nonnegative definite with $\tr C_D=q$, so that $0\leq\theta(D)\leq1$ and
\begin{equation}\label{eq:trace2}
 \tr(C_D^2)=q+q(q-1)\theta(D)^2 .
\end{equation}
The four criteria then become $\det C_D=\det M_D/\det M_F$, $I(D)=1+\tr(C_D^{-1})$, $\eta(D)=\lmin(C_D)$, which equals $\lmin(\mathcal C_D)$ because $\lmin(C_D)\leq1$, and $\Delta(D)$; for a singular design $I(D)=+\infty$ and $\eta(D)=0$. Table~\ref{tab:criteria} lists them with their values for a foldover OofA-OA. \citet{tsai2025} compares OofA-OAs through squared differences between the $J$-characteristics of a design and those of the full design, which for the second moments is the distance between $G_D$ and $\Sigma$ without scaling. The scaling by $\Sigma$ in \eqref{eq:discrepancy} ties the distance to the eigenvalues of $C_D$, and through them to the other criteria.

\begin{table}[tb]
\centering
\caption{Criteria for a foldover design $D$. The last column gives the value for a foldover OofA-OA.}
\label{tab:criteria}
\begin{tabular}{@{}llcc@{}}
\toprule
Criterion & Expression & Optimise & OofA-OA\\
\midrule
Relative determinant & $\det C_D=\Deff(D)^p$ & Maximise & $1$\\
Prediction & $I(D)=1+\tr(C_D^{-1})$ & Minimise & $p$\\
Minimum relative efficiency & $\eta(D)=\lmin(C_D)$ & Maximise & $1$\\
Squared distance & $\Delta(D)=\tr(C_D^2)-q$ & Minimise & $0$\\
\bottomrule
\end{tabular}
\end{table}

The distance can be computed without $\Sigma^{-1/2}$. From \eqref{eq:full} and $(B^\top B)^2=mB^\top B$,
\[
\Sigma^{-1}=3\left(I_q-\frac{B^\top B}{m+1}\right),
\]
and for two orders $a$ and $b$ we define
\begin{equation}\label{eq:kernel}
 \kappa(a,b)=z(a)^\top\Sigma^{-1}z(b)
 =q-6K(a,b)+\frac{6S(a,b)}{m+1},
\end{equation}
where $K(a,b)$ is the Kendall tau distance, the number of pairs of components on whose order $a$ and $b$ disagree, and $S(a,b)=\sum_j\{\pos_a(j)-\pos_b(j)\}^2$. The second equality holds because
\[
 z(a)^\top z(b)=q-2K(a,b),\qquad s(a)^\top s(b)=\frac{m(m^2-1)}3-2S(a,b),
\]
where $s(a)=Bz(a)$ has $j$th entry $m+1-2\pos_a(j)$. In particular $\kappa(a,a)=q$ for every order, so that $x(a)^\top M_F^{-1}x(a)=p$, as in Lemma~1 of \citet{tsai2022}, and $\tr\mathcal C_D=p$ for every design. For a foldover design $\tr C_D=q$ and
\[
 \Delta(D)=\frac1{h^2}\sum_{a,b\in H}\kappa(a,b)^2-q,
\]
where the sum is over all ordered pairs of orders of $H$, including $a=b$; for a design measure the terms are weighted by $w_aw_b$. \citet{shao2026} base space-filling criteria for foldover designs on the Kendall tau distance alone. Both distances enter $\kappa$ because of the scaling by $\Sigma$.

\section{Adding runs to an orthogonal array}\label{sec:adding}

Our main result on run sizes next to those of an OofA-OA is the following.

\begin{theorem}\label{thm:mainaugmentation}
Let $m\geq3$ and let $D_0$ be an OofA-OA with $n_0$ runs.
\begin{enumerate}
\item[(i)] For every order $a$, the design $D_0\cup\{a\}$ is $D$-optimal among all designs with $n_0+1$ runs under the PWO model.
\item[(ii)] If $D_0$ is a foldover design, then for every order $a$ the design $D_0\cup\{a,\rev a\}$ is $D$-optimal among all foldover designs with $n_0+2$ runs under the PWO model.
\end{enumerate}
In both cases the design also minimises $I(D)$ and $\tr\{(\mathcal C_D-I_p)^2\}$ and maximises $\lmin(\mathcal C_D)$ in the same class of designs. The competing designs need not contain $D_0$ and may contain repeated orders.
\end{theorem}

The proof is completed in Section~\ref{sec:fixedsize}, where the optimal values of the four criteria are also given. For $m\geq4$, part (i) settles every run size $n\equiv1\pmod{12}$ at which an OofA-OA with $n-1$ runs exists, and part (ii) the foldover designs with $n\equiv2\pmod{24}$ runs when a foldover OofA-OA with $n-2$ runs exists. For three components the OofA-OAs are the multiples of the full design, and the theorem gives an optimal design at every $n\equiv1\pmod6$ and an optimal foldover design whenever the number $h$ of foldover pairs satisfies $h\equiv1\pmod3$.

All additions of one run are equivalent under the four criteria, but they are not equivalent for every purpose. By the Sherman--Morrison formula, when the order $a$ is added to an OofA-OA with $n-1$ runs, the variance of the fitted mean at an order $b$ is proportional to $p-\{1+\kappa(a,b)\}^2/(n+p-1)$. This is smallest at $b=a$ and decreases with $|1+\kappa(a,b)|$. The order of greatest interest, for example a standard protocol, can therefore be added at no cost in terms of the criteria. The added run also supplies a degree of freedom for error, which is welcome for the smallest arrays: a 12-run OofA-OA with five components leaves one residual degree of freedom, and the 13-run design leaves two.
To illustrate, the twelve orders
\begin{equation}\label{eq:base5}
\begin{array}{llllll}
\mathtt{12354}&\mathtt{14253}&\mathtt{14352}&\mathtt{15243}&\mathtt{15342}&\mathtt{21345}\\
\mathtt{24153}&\mathtt{25143}&\mathtt{31245}&\mathtt{32154}&\mathtt{41325}&\mathtt{42315}
\end{array}
\end{equation}
and their reverses form a 24-run foldover OofA-OA with five components, for which $Z_H^\top Z_H=12\Sigma$. Adding the order $12345$ gives a $D$-optimal 25-run design with
\[
 \Deff=\{(24/25)^{10}(35/25)\}^{1/11}\approx0.993499 .
\]
Adding $12345$ and its reverse gives a 26-run design that is optimal among foldover designs. It has $h=13$, the eigenvalues of $C_D$ are $22/13$ once and $12/13$ nine times, and
\[
 \frac{\det M_D}{\det M_F}=\left(\frac{12}{13}\right)^9\frac{22}{13},\qquad
 \Deff(D)\approx 0.982492,\qquad I(D)\approx 11.3409,
\]
so that $I(D)-1=\tr(C_D^{-1})\approx10.3409$. No 26-run foldover design with five components has a smaller average variance of prediction. With the 20-run result of Theorem~\ref{thm:global5}, optimal foldover designs with five components are now known for $n=20$, $24$ and $26$ and for every $n=2(12k+1)$. In the same way, the 48-run foldover OofA-OAs with six and seven components listed in the Supplementary Material give optimal 50-run foldover designs with $\Deff\approx 0.992105$ and $0.989660$, and adding an order to the 24-run OofA-OA with seven components listed there gives a $D$-optimal 25-run design with $\Deff\approx 0.988813$.
\begin{table}[tb]
\centering
\caption{Optimal designs obtained from an OofA-OA with $n_0$ runs. Adding one order gives the design $D_{\rm O}$ with $n_{\rm O}=n_0+1$ runs, and adding the reverses of all its runs gives the foldover design $D_{\rm F}$ with $n_{\rm F}=2(n_0+1)$ runs. All runs of the listed designs are distinct.}
\label{tab:largerexamples}
\begin{tabular}{rrrrrrr}
\toprule
$m$ & $p$ & $n_0$ & $n_{\rm O}$ & $\Deff(D_{\rm O})$
    & $n_{\rm F}$ & $\Deff(D_{\rm F})$\\
\midrule
7  & 22 & 24  & 25  & 0.988813 & 50  & 0.989660\\
8  & 29 & 48  & 49  & 0.995687 & 98  & 0.995946\\
9  & 37 & 48  & 49  & 0.994839 & 98  & 0.995075\\
10 & 46 & 240 & 241 & 0.999654 & 482 & 0.999668\\
\bottomrule
\end{tabular}
\end{table}

Table~\ref{tab:largerexamples} gives further designs of this kind. They start from the 24-run OofA-OA with seven components of \citet[Example~5]{tsai2022}, from the array with eight components and 48 runs in the \texttt{OofAExp} catalogue \citep{tsai2026package}, which attributes it to \citet{tsai2022}, from the array with nine components and 48 runs of \citet{tsai2025}, and from an array with ten components and 240 runs. The last array is the union of two orbits of the symmetric group on five symbols, acting on the ten components through the ten unordered pairs of symbols. The arrays are listed in Section~S1 of the Supplementary Material. \citet{liu2026} prove universal optimality of component orthogonal arrays under a general class of symmetric models, again at the run sizes where such arrays exist.

\section{Lower bounds and the proof of Theorem~\ref{thm:mainaugmentation}}\label{sec:bound}

In Sections~\ref{sec:spectral} to~\ref{sec:allresidues}, $D$ is a foldover design $H\cup\rev H$ or a design measure that gives equal weights to $a$ and $\rev a$.

\subsection{Eigenvalue bounds}\label{sec:spectral}

The following inequalities for a nonnegative definite matrix with given trace and given trace of its square are classical. The determinant bounds are Theorem~1 of \citet{rodin2017}, the bound on the trace of the inverse follows from Lemma~3.1 of \citet{morgansrivastav2000} within the type~1 theory of \citet[Theorem~2.3]{cheng1978}, and the bound on the smallest eigenvalue is Theorem~2.1 of \citet{wolkowiczstyan1980}.

\begin{lemma}\label{lem:spectral}
Let $C$ be a $q\times q$ nonnegative definite matrix with $\tr C=q$, and define $\theta\in[0,1]$ by $\tr(C^2)=q+q(q-1)\theta^2$.
\begin{enumerate}
\item[(i)] $\det C\leq(1-\theta)^{q-1}\{1+(q-1)\theta\}$.
\item[(ii)] $\tr(C^{-1})\geq(q-1)/(1-\theta)+1/\{1+(q-1)\theta\}$ when $\theta<1$ and $C$ is nonsingular.
\item[(iii)] $\lmin(C)\leq1-\theta$.
\item[(iv)] If $0\leq\theta\leq1/(q-1)$, then $\det C\geq(1+\theta)^{q-1}\{1-(q-1)\theta\}$.
\end{enumerate}
The right-hand sides of (i) and (iii) are decreasing in $\theta$ and that of (ii) is increasing, strictly for $0<\theta<1$. For $0<\theta<1$, equality in (i), (ii) or (iii) holds precisely when the eigenvalues of $C$ are $1+(q-1)\theta$ once and $1-\theta$ with multiplicity $q-1$, and for $0<\theta\leq1/(q-1)$ equality in (iv) holds precisely when they are $1-(q-1)\theta$ once and $1+\theta$ otherwise.
\end{lemma}

By~\eqref{eq:trace2}, the lemma applies to $C=C_D$ with $\theta=\theta(D)$. For PWO designs more can be said about the cases of equality, because the eigenvector of the exceptional eigenvalue has to come from an order.

\begin{theorem}\label{thm:equality}
Let $D$ be a foldover design or a design measure that gives equal weights to $a$ and $\rev a$.
\begin{enumerate}
\item[(i)] For $0<\theta<1$, equality holds in Lemma~\ref{lem:spectral}(i), or in (ii), or in (iii), if and only if there is an order $a\in\Omega_m$ such that
\begin{equation}\label{eq:upperequality}
 G_D=(1-\theta)\Sigma+\theta z(a)z(a)^\top .
\end{equation}
\item[(ii)] For $0<\theta\leq1/(q-1)$, equality holds in Lemma~\ref{lem:spectral}(iv) if and only if $G_D=(1+\theta)\Sigma-\theta z(a)z(a)^\top$ for some order $a$.
\end{enumerate}
Every $\theta\in[0,1]$ is attained by a design measure of the form~\eqref{eq:upperequality}.
\end{theorem}

The proof is given in Section~S11 of the Supplementary Material. For an $n\times p$ matrix $X$ of two-level columns that can be chosen freely, where the reference matrix is the identity, the bound (i) at $\theta=1/n$ is the familiar inequality $\det(X^\top X)\leq(n-1)^{p-1}(n+p-1)$ behind the results of \citet{ehlich1964} and \citet{cheng1980}. Theorem~\ref{thm:equality} says which of the matrices attaining the bounds can come from orders. An optimality result at a given run size still needs a lower bound on $\theta(D)$, and this is supplied in Section~\ref{sec:triples}.

Both cases of equality occur. If $D_0$ is a foldover OofA-OA with $h_0$ foldover pairs, adding one pair gives $\theta=1/(h_0+1)$ and
\begin{equation}\label{eq:addone}
 \frac{\det M_{D_+}}{\det M_F}
 =\left(\frac{h_0}{h_0+1}\right)^{q-1}\frac{h_0+q}{h_0+1},
\end{equation}
which attains Lemma~\ref{lem:spectral}(i), whereas deleting one of its pairs, when $h_0\geq q$, gives $\theta=1/(h_0-1)$ and
\begin{equation}\label{eq:deleteone}
 \frac{\det M_{D_-}}{\det M_F}
 =\left(\frac{h_0}{h_0-1}\right)^{q-1}\frac{h_0-q}{h_0-1},
\end{equation}
which attains Lemma~\ref{lem:spectral}(iv). Adding a pair thus reaches the upper bound on the determinant at its value of $\theta$, and deleting a pair reaches the lower bound. We return to this asymmetry in Section~\ref{sec:deletion}.

\subsection{A representation by triples}\label{sec:triples}

For two-level designs the lower bound comes from the parity of the entries of $X^\top X$. The entries of $Z_H^\top Z_H$ do not keep such properties under the scaling by $\Sigma$. The next lemma replaces the PWO factors by signs attached to triples of labels, and the scaled distance becomes a sum of squares to which parity and divisibility apply.

\begin{lemma}\label{lem:triples}
Introduce the label $0$ and let $\mathcal T$ be the set of triples of $\{0,1,\ldots,m\}$, so that $|\mathcal T|=\binom{m+1}3$. For an order $a$ define
\[
 u_{0ij}(a)=z_{ij}(a),\qquad
 u_{ijk}(a)=z_{ij}(a)-z_{ik}(a)+z_{jk}(a)\quad(0<i<j<k),
\]
and write $u(a)=Lz(a)$, where $L$ is the $|\mathcal T|\times q$ matrix defined by these expressions. Then every $u_T(a)$ is $1$ or $-1$, and $u_T(a)$ is the orientation of the triple $T$ in the circular order $(0,a_1,\ldots,a_m)$: it equals $1$ exactly when, starting from the smallest label of $T$ and moving around the circle in the direction of addition, the other two labels of $T$ are met in increasing order. Moreover
\begin{equation}\label{eq:tripleframe}
 L^\top L=(m+1)I_q-B^\top B,\qquad
 \Sigma^{-1}=\frac3{m+1}L^\top L .
\end{equation}
For a collection $H$ of $h$ orders, repetitions allowed, let $A=Z_H^\top Z_H$, let $U=Z_HL^\top$ be the $h\times|\mathcal T|$ matrix of the orientations, and let $R=L\Sigma L^\top$. Then
\begin{equation}\label{eq:triplesquares}
 V(H):=\tr\!\left[\{\Sigma^{-1}(A-h\Sigma)\}^{2}\right]
 =\frac9{(m+1)^2}\|U^\top U-hR\|_F^2 ,
\end{equation}
where $\|\cdot\|_F$ is the Frobenius norm, and $V(H)=h^2\Delta(D)=h^2q(q-1)\theta(D)^2$ for $D=H\cup\rev H$. The matrix $R$ has unit diagonal, $R_{ST}=0$ when the triples $S\neq T$ share at most one label, and $|R_{ST}|=1/3$ when they share two labels.
\end{lemma}
\begin{proof}
The identity $z_{ij}-z_{ik}+z_{jk}=z_{ij}z_{ik}z_{jk}$ of \citet{mee2020} gives $u_T(a)=\pm1$, and comparing with the positions of $i$, $j$, $k$ on the circle through $0$ gives the statement on orientations. Counting the diagonal and off-diagonal entries of $L^\top L$ gives the first identity in~\eqref{eq:tripleframe}, and the second follows from $\Sigma^{-1}=3\{I_q-B^\top B/(m+1)\}$. Substituting~\eqref{eq:tripleframe} into $V(H)$ and using the cyclic invariance of the trace gives~\eqref{eq:triplesquares}. The uniform distribution on orders induces the uniform distribution on circular orders of $\{0,1,\ldots,m\}$, so $R_{ST}=E(u_Su_T)$. If $S\neq T$ share at most one label, exchanging two labels of $S\setminus T$ reverses $u_S$ and preserves $u_T$; the exchange may involve the label $0$, because the symmetric group on $\{0,1,\ldots,m\}$ acts on circular orders and preserves their uniform distribution. Hence $R_{ST}=0$. If they share two labels, they lie in a four-element set $Q$, treated in the proof of Lemma~\ref{lem:integerdiscrepancy}, which gives $|R_{ST}|=1/3$.
\end{proof}

\begin{lemma}\label{lem:integerdiscrepancy}
Let $H$ be a collection of $h$ orders on $m\geq3$ components, repetitions allowed, with $h$ odd and not divisible by three, and let $V(H)$ be as in~\eqref{eq:triplesquares}. Then $V(H)\geq q(q-1)$; equivalently, the foldover design $D=H\cup\rev H$ has $h\theta(D)\geq1$. For $m=3$ the conclusion holds whenever $h$ is not divisible by three, whether $h$ is odd or even.
\end{lemma}
\begin{proof}
The diagonal entries of $U^\top U-hR$ vanish, so only the $|\mathcal T|(|\mathcal T|-1)$ ordered pairs of distinct triples contribute to~\eqref{eq:triplesquares}. If $S,T$ share at most one label, the entry is $\sum_{a\in H}u_S(a)u_T(a)$, a sum of $h$ signs; since $h$ is odd it is an odd integer and its square is at least one.

The pairs sharing two labels are grouped by four-element sets. On $Q=\{i,j,k,\ell\}$ with $i<j<k<\ell$, put $w=(u_{jk\ell},-u_{ik\ell},u_{ij\ell},-u_{ijk})^\top$; the defining expressions give $\mathbf1^\top w=0$, so $w$ has two positive and two negative entries, and $w_1w_2=w_3w_4$, $w_1w_3=w_2w_4$, $w_1w_4=w_2w_3$ and $w_1w_2+w_1w_3+w_1w_4=-1$. A circular order of four labels is determined by the orientations of its four triples, so the six circular orders on $Q$ give the six balanced sign patterns of $w$ once each, and $E(w_rw_s)=-1/3$ for $r\neq s$, which is the last statement of Lemma~\ref{lem:triples}. Let $b_1,b_2,b_3$ be the sums over $H$ of the three products $w_1w_2$, $w_1w_3$, $w_1w_4$. Each is congruent to $h$ modulo two, their sum is $-h$, and the six unordered off-diagonal positions indexed by $Q$ have values $\pm(b_r+h/3)$, each value occurring twice. Since three does not divide $h$,
\begin{equation}\label{eq:localrounding}
 \sum_{r=1}^3(b_r+h/3)^2\geq\frac83 ,
\end{equation}
because the three integers have a common parity and so cannot differ by one, the sum is minimised when they differ by at most two, otherwise moving two units from the largest to the smallest reduces it, and the deviations from $-h/3$ are then $2/3,2/3,-4/3$ up to sign and permutation. Each of the $\binom{m+1}4$ sets $Q$ thus contributes at least $4\times8/3=32/3$ to the squared norm, from its twelve ordered pairs of triples.

There are $12\binom{m+1}4$ ordered pairs sharing two labels and $|\mathcal T|(|\mathcal T|-1)-12\binom{m+1}4$ sharing at most one. Both contributions apply, so
\[
 \|U^\top U-hR\|_F^2
 \geq |\mathcal T|(|\mathcal T|-1)-12\binom{m+1}{4}+\frac{32}{3}\binom{m+1}{4}
 =\frac{(m+1)^2}{9}q(q-1),
\]
the last equality being a polynomial identity in $m$. Multiplying by $9/(m+1)^2$ gives $V(H)\geq q(q-1)$, and the statement on $\theta(D)$ follows from $V(H)=h^2q(q-1)\theta(D)^2$. Repeated orders, and a foldover pair of which both orders are in $H$, do not affect any step, because every entry is a sum over the rows of $H$ and $z(\rev a)=-z(a)$. For $m=3$ every two of the four triples on $\{0,1,2,3\}$ share two labels, so the parity contribution is void and the four-element set alone gives $\|U^\top U-hR\|_F^2\geq32/3$, that is $V(H)\geq6=q(q-1)$, for every $h$ not divisible by three.
\end{proof}

Equality in Lemma~\ref{lem:integerdiscrepancy} holds exactly when every term attains its bound: every sum $\sum_{a\in H}u_S(a)u_T(a)$ for triples sharing at most one label has absolute value one, and on every four-element set the deviations $b_r+h/3$ in~\eqref{eq:localrounding} are $2/3,2/3,-4/3$ up to a common sign and permutation. When $h$ is even only the four-element sets contribute, and when $h$ is divisible by three only the parity terms do. The weaker bounds obtained in these cases are given in Section~S3 of the Supplementary Material, and they are improved next.

\subsection{Bounds for every number of foldover pairs}\label{sec:allresidues}

Bounds for all $h$ follow once the case of four components is solved completely. For $r\in\{0,\ldots,11\}$ define
\begin{equation}\label{eq:residuevalues}
 (\mu_0,\ldots,\mu_{11})=
 \left(0,30,48,\frac{1566}{25},\frac{1632}{25},\frac{294}{5},
 \frac{216}{5},\frac{294}{5},\frac{1632}{25},\frac{1566}{25},48,30\right).
\end{equation}

\begin{theorem}\label{thm:m4discrepancyall}
For four components and every integer $h\geq0$, the minimum of $V(H)$ over all half-designs $H$ with $h$ orders is $\mu_r$, where $r\equiv h\pmod{12}$ and $0\leq r<12$. Write $h=12k+r$. Every foldover design attaining the minimum consists of $k$ copies of each of the twelve foldover pairs and $r$ further distinct pairs. For $r=1$ the further pair is arbitrary; for $r=2$ the two pairs form one of six couples listed in Section~S9 of the Supplementary Material; for $r=3,4,5,6$ their orders are all even permutations or all odd permutations; and for $r=7,\ldots,11$ they are the pairs left out by one of the sets just described for $12-r$. In particular, the numbers of foldover designs attaining the minimum for $r=0,\ldots,11$ are
\[
 1,12,6,40,30,12,2,12,30,40,6,12,
\]
independently of $k$.
\end{theorem}

For $m\geq4$ put $\nu=m+1$ and
\begin{equation}\label{eq:residueboundconstant}
\begin{split}
 B_m(h)={}&\epsilon_h q(q-1)
 +\frac{5\binom{\nu}{4}}{\nu^2}
        \{\mu_{h\bmod12}-30\epsilon_h\}\\
 &+\frac{72(\nu-5)\binom{\nu}{4}}{\nu^2}
       \mathbf1_{\{h\equiv2\ (\mathrm{mod}\ 4)\}},
 \qquad \epsilon_h=\mathbf1_{\{h\text{ odd}\}}.
\end{split}
\end{equation}

\begin{theorem}\label{thm:allresiduediscrepancy}
For every $m\geq4$ and every collection $H$ of $h$ orders, repetitions allowed,
\[
 V(H)\geq B_m(h).
\]
For $m=4$ the bound is attained for every $h$. For every $m\geq4$, $B_m(h)=0$ exactly when $h\equiv0\pmod{12}$, and $B_m(h)=q(q-1)$ when $h\equiv1,11\pmod{12}$. For $h\equiv2,\ldots,10\pmod{12}$ the bound is strictly larger than the bounds obtained from parity and from the four-element sets alone.
\end{theorem}

\begin{theorem}\label{thm:oddminimumstructure}
Let $m\geq4$ and let $h$ be odd. Then every collection $H$ of $h$ orders, repetitions allowed, has $V(H)\geq q(q-1)$. Equality holds if and only if one of the following two alternatives holds, for some order $a$:
\begin{align*}
 h&\equiv1\pmod{12}, &
 Z_H^\top Z_H&=(h-1)\Sigma+z(a)z(a)^\top;\\
 h&\equiv11\pmod{12}, &
 Z_H^\top Z_H&=(h+1)\Sigma-z(a)z(a)^\top.
\end{align*}
These conditions concern the matrix $Z_H^\top Z_H$ and do not say that $H$ contains the order $a$ or an OofA-OA.
\end{theorem}

The three theorems are proved in Section~S9 of the Supplementary Material. For Theorem~\ref{thm:m4discrepancyall}, $25V(H)$ is written as a quadratic form in the multiplicities of the twelve foldover pairs and minimised over the integers. Theorem~\ref{thm:allresiduediscrepancy} follows by applying Theorem~\ref{thm:m4discrepancyall} to the circular orders induced on every five labels of $\{0,1,\ldots,m\}$, with one more parity condition on six labels when $h\equiv2\pmod4$. Theorem~\ref{thm:oddminimumstructure} then shows that the bound $V(H)\geq q(q-1)$ of Lemma~\ref{lem:integerdiscrepancy} holds for every odd $h$ when $m\geq4$, and that it is strict unless $h\equiv\pm1\pmod{12}$. A matrix attaining it at $h\equiv1\pmod{12}$ has one larger eigenvalue on the scale of $\Sigma$ and $q-1$ equal smaller ones, and attains the upper bound of Lemma~\ref{lem:spectral}(i); at $h\equiv11\pmod{12}$ it has one smaller eigenvalue and $q-1$ equal larger ones, and attains the lower bound of Lemma~\ref{lem:spectral}(iv). Whether a design attains $q(q-1)$ at a given $h$ is a separate question.

Through Lemma~\ref{lem:spectral} these results bound the criteria at every run size. For a foldover design with $h>0$ pairs put $t_{m,h}=h^{-1}\{B_m(h)/q(q-1)\}^{1/2}$. The right-hand sides of Lemma~\ref{lem:spectral}(i)--(iii) with $t_{m,h}$ in place of $\theta$ then bound $\det C_D$, $\tr(C_D^{-1})$ and $\lmin(C_D)$, and corresponding bounds hold for an arbitrary design with $n$ runs, with $p$ in place of $q$ and $B_m(n)+2q\epsilon_n$ in place of $B_m(h)$ (Corollary~S3 in the Supplementary Material). For $m>4$ these bounds need not be attained, and for small $h$ they can be far from the optimum. With five components and ten pairs the bound on the efficiency is $0.963$, whereas the optimal foldover design of Theorem~\ref{thm:global5} has efficiency $0.898$.

\subsection{Proof of Theorem~\ref{thm:mainaugmentation}}\label{sec:fixedsize}

Combining the lower bounds on $V(H)$ with Lemma~\ref{lem:spectral} gives the result for foldover designs.

\begin{theorem}\label{thm:fixedsize}
Let $m\geq3$, and let $D$ be a foldover design with $h>1$ foldover pairs, counted with multiplicity, where $h$ is odd when $m\geq4$ and $h$ is not divisible by three when $m=3$. Then
\begin{equation}\label{eq:fixedsizebound}
\begin{gathered}
 \frac{\det M_D}{\det M_F}\leq\left(\frac{h-1}{h}\right)^{q-1}\frac{h+q-1}{h},\qquad
 \tr(C_D^{-1})\geq\frac{(q-1)h}{h-1}+\frac h{h+q-1},\\
 \lmin(C_D)\leq\frac{h-1}h,\qquad
 \tr(C_D^2)\geq q+\frac{q(q-1)}{h^2},
\end{gathered}
\end{equation}
where $\tr(C_D^{-1})$ is taken as $+\infty$ for a singular design. Equality holds in any of the first three bounds if and only if, for an order $a$,
\begin{equation}\label{eq:fixedsizeequality}
 Z_H^\top Z_H=(h-1)\Sigma+z(a)z(a)^\top ,
\end{equation}
which requires $h\equiv1\pmod{12}$ when $m\geq4$ and $h\equiv1\pmod3$ when $m=3$. Equality in the fourth holds if and only if $V(H)=q(q-1)$. In particular, adding any foldover pair to a foldover OofA-OA gives a design that attains all four bounds.
\end{theorem}
\begin{proof}
Theorem~\ref{thm:oddminimumstructure} for $m\geq4$ and Lemma~\ref{lem:integerdiscrepancy} for $m=3$ give $V(H)\geq q(q-1)$, that is $\theta\geq1/h$, and the right-hand sides in Lemma~\ref{lem:spectral}(i)--(iii) are monotone in $\theta$, so the first three bounds follow by evaluating them at $\theta=1/h$; the fourth is~\eqref{eq:trace2}. Equality in any of the first three requires $\theta=1/h$ together with equality in Theorem~\ref{thm:equality}(i), which is~\eqref{eq:fixedsizeequality}; conversely, a design satisfying~\eqref{eq:fixedsizeequality} has $\theta=1/h$ and the extremal eigenvalues, so it attains the three bounds. A foldover OofA-OA with $h-1$ pairs has $Z^\top Z=(h-1)\Sigma$, so any added pair gives~\eqref{eq:fixedsizeequality}; its number of pairs $h-1$ is a multiple of twelve when $m\geq4$ by Corollary~\ref{cor:sizes} and a multiple of three when $m=3$, so $h$ satisfies the hypotheses. For $m\geq4$, Proposition~\ref{prop:counts} applied to~\eqref{eq:fixedsizeequality} with $\varepsilon=1$ and $t=1$ gives $h-1\in12\mathbb Z$. For $m=3$ there are three foldover pairs, and the off-diagonal entries of~\eqref{eq:fixedsizeequality} force the multiplicities of the pairs other than that of $a$ to equal $(h-1)/3$, so $h\equiv1\pmod3$.
\end{proof}

Theorem~\ref{thm:fixedsize} contains Theorem~\ref{thm:mainaugmentation}(ii). Repeating a foldover OofA-OA with $h_0$ pairs $k$ times and adding any pair also attains the four bounds, with $2(kh_0+1)$ runs. For $m\geq4$ and odd $h\not\equiv1\pmod{12}$ the first three bounds are not attained and need not be close to the optimum. With four components, a complete enumeration of the foldover designs, repeated pairs allowed, gives a largest $\Deff\approx0.9077$ at $h=7$ and $0.9759$ at $h=11$, against the bounds $0.9674$ and $0.9856$, and each optimal design is unique up to relabelling. The first consists of the even or the odd permutations and one further pair, and the second is the full design with one pair deleted (Table~S1 in the Supplementary Material).

For a design with an odd number of runs the column sums of $Z_D$ cannot all vanish, and the same representation bounds them as well. This extends the comparison from foldover designs to all designs.

\begin{corollary}\label{cor:ordinaryfixed}
Let $m\geq3$ and let $D$ be any design with $n$ runs, where $n>1$ is odd and, when $m=3$, not divisible by three. Then
\[
\begin{gathered}
 \frac{\det M_D}{\det M_F}\leq\left(\frac{n-1}{n}\right)^{p-1}\frac{n+p-1}{n},\qquad
 \tr(M_D^{-1}M_F)\geq\frac{(p-1)n}{n-1}+\frac n{n+p-1},\\
 \lmin(\mathcal C_D)\leq\frac{n-1}n,\qquad
 \tr(\mathcal C_D^2)\geq p+\frac{p(p-1)}{n^2},
\end{gathered}
\]
where $\tr(M_D^{-1}M_F)=I(D)$ is taken as $+\infty$ for a singular design. Equality holds in any of the first three bounds if and only if
\[
 X_D^\top X_D=(n-1)M_F+x(a)x(a)^\top\quad\text{for an order }a,
\]
which requires $n\equiv1\pmod{12}$ when $m\geq4$ and $n\equiv1\pmod6$ when $m=3$. Whenever an OofA-OA with $n-1$ runs exists, adding any one order to it attains all four bounds.
\end{corollary}

The proof is given in Section~S11 of the Supplementary Material. With $A=Z_D^\top Z_D$ and $b_D=Z_D^\top\mathbf1$ it rests on the identity
\[
 n^2\tr\{(\mathcal C_D-I_p)^2\}
 =\tr[\{\Sigma^{-1}(A-n\Sigma)\}^2]+2b_D^\top\Sigma^{-1}b_D.
\]
The first term is at least $q(q-1)$ by Theorem~\ref{thm:oddminimumstructure} or Lemma~\ref{lem:integerdiscrepancy}. The second is at least $2q$, because by~\eqref{eq:tripleframe} it equals $6\|Lb_D\|^2/(m+1)$ and every entry of $Lb_D$ is a sum of $n$ signs, hence an odd integer. Lemma~\ref{lem:spectral} is then applied to $\mathcal C_D$ with $p$ in place of $q$. Corollary~\ref{cor:ordinaryfixed} contains Theorem~\ref{thm:mainaugmentation}(i), and the proof of Theorem~\ref{thm:mainaugmentation} is complete.

The forms of the matrices attaining the bounds can be recognised from counts on four components. On a set $Q\subseteq\{1,\ldots,m\}$ of four components the 24 relative orders form twelve foldover pairs. For an order $a$, let $[a|_Q]$ be the foldover pair containing the relative order of the components of $Q$ in $a$, and let $c_{Q,b}$ be the number of orders of $H$ whose relative order on $Q$ belongs to the foldover pair $b$. These counts do not change when an order of $H$ is replaced by its reverse, and they add up to $h$.

\begin{proposition}\label{prop:counts}
For $m\geq4$, let $D$ be a foldover design with $h$ foldover pairs, counted with multiplicity, and fix an order $a$, a number $\theta>0$ and $\varepsilon\in\{1,-1\}$. Put $t=h\theta$. Then
\begin{equation}\label{eq:signedequality}
 G_D=(1-\varepsilon\theta)\Sigma+\varepsilon\theta z(a)z(a)^\top
\end{equation}
if and only if, for every set $Q$ of four components and every foldover pair $b$ of relative orders on it,
\begin{equation}\label{eq:countform}
 c_{Q,b}=\frac{h-\varepsilon t}{12}
 +\varepsilon t\,\mathbf1_{\{b=[a|_Q]\}}.
\end{equation}
In particular, \eqref{eq:signedequality} requires $t\in\mathbb Z$ and $h-\varepsilon t\in12\mathbb Z$, and for $\varepsilon=-1$ also $h\geq11t$.
\end{proposition}

The proof is given in Section~S11 of the Supplementary Material. As a check, in the 98-run foldover design with eight components of Table~\ref{tab:largerexamples}, for which the added order is $12345678$, on every set of four components the foldover pair of the increasing order occurs five times in $H$ and each of the other eleven pairs four times. The case $\theta=0$ is known. In a foldover design all moments of odd order of the PWO factors vanish, so a foldover OofA-OA has strength three \citep[Corollary~3]{tsai2025}, and the results of \citet[Lemma~1, Remark~1 and Theorem~2]{zhaodong2022} on arrays of strength three give equal frequencies of the relative orders on four components and a run size divisible by 24. Every foldover pair on four components then occurs $h/12$ times in $H$. Together with the requirement $h\geq q$ for estimability this gives the following bounds.

\begin{corollary}\label{cor:sizes}
A foldover OofA-OA with $m\geq4$ components has a run size divisible by 24, and hence
\[
 n\geq24\left\lceil\frac{m(m-1)}{24}\right\rceil .
\]
These lower bounds are 24, 48, 48, 72, 72 and 96 for $m=5,\ldots,10$, and are attained for $m=5,6,7$. For OofA-OAs that need not be foldover designs, the corresponding bounds $12\lceil p/12\rceil$ are 12, 24, 24, 36, 48 and 48.
\end{corollary}

The array in~\eqref{eq:base5} and the arrays with six and seven components in Section~S1 of the Supplementary Material attain the first three bounds. Thus no foldover OofA-OA with six components has 36 runs, although 36 is a multiple of 12. The half-design listed for the 48-run array with six components is not itself an OofA-OA, but replacing some of its orders by their reverses gives a 24-run OofA-OA, and the 48-run array consists of that OofA-OA and the reverses of its runs. Whether foldover OofA-OAs exist with 72 runs for eight or nine components and with 96 runs for ten is open. Proposition~\ref{prop:counts} also explains the congruences in Theorems~\ref{thm:oddminimumstructure} and~\ref{thm:fixedsize}. With $t=1$, the form with $\varepsilon=1$ requires $h\equiv1\pmod{12}$ and the form with $\varepsilon=-1$ requires $h\equiv-1\pmod{12}$, and on each set of four components one foldover pair then occurs once more, or once less, than each of the other eleven. This is why the first three bounds of Theorem~\ref{thm:fixedsize} cannot be attained for odd $h\not\equiv1\pmod{12}$, whereas the fourth is attained at $h\equiv11\pmod{12}$ by a foldover OofA-OA with one pair deleted. In the same way, the fourth bound of Corollary~\ref{cor:ordinaryfixed} is attained at $n\equiv11\pmod{12}$ by an OofA-OA with one run deleted.

Equality in Theorem~\ref{thm:fixedsize} and Corollary~\ref{cor:ordinaryfixed} is a statement about the moment matrix, and a design attaining the bounds need not contain the order $a$ when orders are repeated. For four components and 13 runs a complete enumeration finds 24 designs that attain the bounds of Corollary~\ref{cor:ordinaryfixed} and contain neither the order $a$ of their moment matrix nor any OofA-OA, which parallels an observation of \citet{hameed2025} for two-level designs. In contrast, every foldover design with four components and 13 pairs that attains the bounds of Theorem~\ref{thm:fixedsize} is the full design with one added pair. These enumerations, the smallest run size at which the bounds are attained for nine components, and a further family of designs attaining them, obtained by ranking the rows of orthogonal arrays of strength four, are given in Sections~S7 and~S6 of the Supplementary Material.
\section{Deleting a run or a foldover pair}\label{sec:deletion}

One might expect that deleting a run from an OofA-OA is as good as adding one. This is not so in general. Let $D_0$ be an OofA-OA with $n_0>p$ runs, and let $D_-$ be the design obtained by deleting one of its runs, the order $a$ say. Then $X_{D_-}^\top X_{D_-}=n_0M_F-x(a)x(a)^\top$, the matrix $\mathcal C_{D_-}$ has the eigenvalue $n_0/(n_0-1)$ with multiplicity $p-1$ and the eigenvalue $(n_0-p)/(n_0-1)$, and
\begin{equation}\label{eq:deleterun}
 \frac{\det M_{D_-}}{\det M_F}
 =\left(\frac{n_0}{n_0-1}\right)^{p-1}\frac{n_0-p}{n_0-1},
\end{equation}
whichever run is deleted. The design $D_-$ attains the fourth bound of Corollary~\ref{cor:ordinaryfixed}, so no design with $n_0-1$ runs is closer to the full design. Its determinant, however, is the lower bound of Lemma~\ref{lem:spectral}(iv), applied to $\mathcal C_D$ with $p$ in place of $q$. Among all nonnegative definite matrices with the same trace and at the same distance from the identity, $\mathcal C_{D_-}$ has the smallest determinant, whereas the matrix of an OofA-OA with an added run has the largest.

\citet[Theorem~2]{tsai2022} states that $D_-$ is $D$-optimal among all designs with $n_0-1$ runs. The proof given there uses the inequality $x(a)^\top(X^\top X)^{-1}x(a)\geq p/n$ for a run $a$ of an arbitrary design with $n$ runs and model matrix $X$. For an OofA-OA there is equality for every run, but in general $p/n$ is only the average of these quantities over the runs of the design. The statement is true in the cases that we could enumerate, namely three components with 5, 11 and 17 runs and four components with 11 runs, repeated runs being allowed, and four components with 23 distinct runs (Section~S7 of the Supplementary Material). It fails in the three cases with a deleted run in Table~\ref{tab:deletioncomparison}.

\begin{table}[tb]
\centering
\caption{An OofA-OA with one run deleted, or a foldover OofA-OA with one foldover pair deleted, compared with another design with the same number $n$ of runs. The last column is the ratio of the determinant of the competing design to that of the design with the deletion.}
\label{tab:deletioncomparison}
\begin{tabular}{@{}rrlccc@{}}
\toprule
 & & & \multicolumn{2}{c}{$\Deff$} & \\
\cmidrule(lr){4-5}
$m$ & $n$ & Deleted & After deletion & Competing design & Ratio\\
\midrule
5 & 11 & One run of 12 & 0.870324 & 0.902665 & $121/81$\\
6 & 23 & One run of 24 & 0.974234 & 0.974621 & $789/784$\\
7 & 23 & One run of 24 & 0.932030 & 0.945984 & $355/256$\\
7 & 46 & One pair of 24 & 0.947532 & 0.955811 & $155/128$\\
\bottomrule
\end{tabular}
\end{table}

With five components $p=11$, so that a 12-run OofA-OA with one run deleted is a saturated design, and~\eqref{eq:deleterun} gives $|\det X_{D_-}|=36\,864$. The eleven orders
\begin{equation}\label{eq:eleven}
\begin{array}{llllll}
\mathtt{12345}&\mathtt{15324}&\mathtt{15423}&\mathtt{24315}&\mathtt{25341}&\mathtt{32514}\\
\mathtt{34215}&\mathtt{35412}&\mathtt{41325}&\mathtt{45213}&\mathtt{45312}&
\end{array}
\end{equation}
have $|\det X_D|=45\,056$, so that $\det(X_D^\top X_D)$ is larger by the factor $(11/9)^2$. This design is also better under the other two criteria that depend on the eigenvalues. It has $I(D)=307/22\approx13.95$ against $121/6\approx20.17$, and the smallest eigenvalue of $\mathcal C_D$ is $4/11$ against $1/11$. It is farther from the full design, with $n^2\tr\{(\mathcal C_D-I_p)^2\}=206$ against the smallest possible value $p(p-1)=110$. For six and seven components the 24-run OofA-OAs of \citet[Examples~3 and~5]{tsai2022} with one run deleted are improved in the same way by the 23-run designs listed in Section~S11.5 of the Supplementary Material. The value $45\,056$ is not new: \citet[Table~6.1]{winker2020} report the efficiency $0.90267$ for five components and 11 runs, and \citet{wangwang2023} list another 11-run design with this determinant. To our knowledge it has not been noted that these values exceed what Theorem~2 of \citet{tsai2022} allows. No exchange of a single run increases the determinant of any of these designs with a run deleted, so an exchange algorithm that accepts only improvements stops there, although for five components a sequence of exchanges that begins with two exchanges of unchanged determinant leads to $|\det X_D|=45\,056$ (Section~S11.5 of the Supplementary Material). The competing designs in Table~\ref{tab:deletioncomparison} were found by computer search from random starting designs, and we do not claim that they are optimal.

The situation for two-level designs is similar. With $n\equiv3\pmod4$ runs and $p$ parameters, an orthogonal array with one run deleted is $D$-optimal when $n>2p-5$ and need not be so when $p$ is larger \citep{galilkiefer1980}. The cases above fit this pattern, and so does the case of five components and 23 runs, where $n>2p-5$ and exchange searches found no design better than $D_-$. Whether the same inequality governs the PWO model, we do not know.

Foldover pairs behave in the same way. Let $D_0$ now be the 48-run foldover OofA-OA with seven components that consists of the 24-run OofA-OA of \citet[Example~5]{tsai2022} and the reverses of its runs, so that $q=21$ and $h_0=24$. By~\eqref{eq:deleteone}, deleting any one of its pairs gives
\begin{equation}\label{eq:deletionbenchmark}
 \det(Z_H^\top Z_H)=2^{78},\qquad
 \Deff(D_-) \approx 0.947532,\qquad \theta=\frac1{23}.
\end{equation}
The 23 orders listed in Section~S11.5 of the Supplementary Material, together with their reverses, form a foldover design with 46 distinct runs and
\begin{equation}\label{eq:counterexample}
 \det(Z_H^\top Z_H)=155\cdot2^{71},\qquad
 \Deff(D) \approx 0.955811,\qquad \theta \approx 0.0475.
\end{equation}
The design was found by computer search, and with a suitable choice of one order from each pair its 23 pairs give the 23-run design with seven components of Table~\ref{tab:deletioncomparison}.

This foldover design is not a foldover OofA-OA with one pair deleted, since $24\Sigma-Z_H^\top Z_H$ has rank three, and it has no foldover pair in common with $D_0$. Nor does any exchange of a single pair increase the determinant of a design obtained from $D_0$ by deleting a pair, when designs with distinct runs are considered (Section~S5 of the Supplementary Material). Whether the design~\eqref{eq:counterexample} is optimal among 46-run foldover designs is open. Here too the smaller distance from the full design goes with the smaller determinant. The design with a pair deleted has the smallest value of $V(H)$ that Lemma~\ref{lem:integerdiscrepancy} allows at $h=23$, namely $V=420=q(q-1)$, whereas the design~\eqref{eq:counterexample} has $V=501$. In fact, by Theorem~\ref{thm:oddminimumstructure} every design with $V(H)=420$ at $h=23$ has $\det(Z_H^\top Z_H)=2^{78}$, so none of them is $D$-optimal.
\section{Foldover designs}\label{sec:consequences}

Theorem~\ref{thm:mainaugmentation}(ii) is an optimality result within the class of foldover designs. In this section we show when that class is the right one, what the restriction costs under the PWO model, and how foldover designs behave away from the run sizes covered by Theorem~\ref{thm:mainaugmentation}.

\subsection{A mean with a reversal-invariant term}\label{sec:exactnuisance}

\citet{mee2020} observed that in a design made of foldover pairs the estimators of the PWO effects are not affected by interactions of two PWO factors. More generally, suppose that the mean contains, in addition to the PWO effects, a term that does not change when the order is reversed,
\begin{equation}\label{eq:nuisance}
 E\{Y(a)\}=z(a)^\top\beta+g(a),\qquad g(\rev a)=g(a),
\end{equation}
where the values of $g$ on different foldover pairs are unknown and unrelated. The intercept and every interaction of an even number of PWO factors are of this form, and the model does not say which of them are present. Eliminating $g$ within the pairs, as block effects are eliminated in a paired comparison, leaves the following information on $\beta$.

\begin{proposition}\label{thm:nuisance}
For each foldover pair in a design, choose one of its orders $a$ and let $n_a^+$ and $n_a^-$ be the numbers of runs at $a$ and at $\rev a$, with $t_a=n_a^++n_a^->0$. Under model~\eqref{eq:nuisance}, the information matrix for $\beta$ after eliminating $g$ is
\begin{equation}\label{eq:nuisanceinfo}
 \mathcal I_\beta=\frac1{\sigma^2}\sum_a
 \frac{4n_a^+n_a^-}{t_a}\,z(a)z(a)^\top,
\end{equation}
so that a pair with $n_a^+=0$ or $n_a^-=0$ contributes nothing. The PWO effects are estimable if and only if the vectors $z(a)$ of the pairs with both orders present span $\R^q$. For given totals $t_a$, the matrix $\mathcal I_\beta$ is largest in the Loewner order when each total is divided as equally as possible. The smallest number of runs with which the PWO effects are estimable is $2q$, and it is attained with distinct runs for every $m\geq3$.
\end{proposition}
\begin{proof}
The Schur complement within pairs gives~\eqref{eq:nuisanceinfo}, and $4n_a^+n_a^-/t_a=t_a-(n_a^+-n_a^-)^2/t_a$ gives the statement on the division of the totals. At least $q$ pairs with both orders present are necessary. A set of $q$ orders that contains no foldover pair and has linearly independent PWO vectors is obtained by reversing $(1,\ldots,m)$ through $q$ interchanges of adjacent components, as shown in Section~S2 of the Supplementary Material.
\end{proof}

When each pair is run once in each order, the half-difference $\{Y(a)-Y(\rev a)\}/2$ has mean $z(a)^\top\beta$ and variance $\sigma^2/2$, and the $h$ half-differences carry the information $2hG_D/\sigma^2$. More generally, a foldover design has $\mathcal I_\beta=nG_D/\sigma^2$ by~\eqref{eq:nuisanceinfo}. Since $\det M_D=\det G_D$ for a foldover design, the $D$-criterion ranks foldover designs in the same way under model~\eqref{eq:nuisance}, with $q$ parameters of interest, and under the PWO model with $p$ parameters. \citet{wangmee2022} use foldover designs to reduce the bias caused by a misspecified model, and~\eqref{eq:nuisanceinfo} expresses this protection for an arbitrary reversal-invariant term and unequal replication. Its limits are discussed in Section~S2 of the Supplementary Material.

The best allocation of $n$ runs under~\eqref{eq:nuisance} does not depend on the special form of the PWO factors, and we state it generally. Suppose that there are finitely many candidate design points, arranged in pairs $\{a,\rho a\}$, where $\rho a\neq a$ and $\rho(\rho a)=a$, that the regression vectors satisfy $f(\rho a)=-f(a)\in\R^q$, are nonzero and span $\R^q$, where $q$ now stands for the dimension of $\beta$, and that
\[
 E\{Y(a)\}=f(a)^\top\beta+g(a),\qquad g(\rho a)=g(a),
\]
where the values of $g$ on different pairs are unrelated fixed effects and the errors are uncorrelated with variance $\sigma^2$. Choose one point $a_i$ from each pair, write $f_i=f(a_i)$, and let $r_i$ and $s_i$ be the numbers of runs at $a_i$ and $\rho a_i$. As in Proposition~\ref{thm:nuisance}, the information matrix for $\beta$ after eliminating $g$ is $\sigma^{-2}J$, where
\begin{equation}\label{eq:generalnuisanceJ}
 J=\sum_i \psi(r_i,s_i)f_if_i^\top,
 \qquad \psi(r,s)=\begin{cases}4rs/(r+s),&r+s>0,\\0,&r=s=0.\end{cases}
\end{equation}
Designs are compared at the same total number of runs $n=\sum_i(r_i+s_i)$. Replication is allowed, and the totals $r_i+s_i$ of the pairs are not fixed in advance. For a given total $t=r+s$ the largest value of $\psi(r,s)$ is
\[
 \omega(t)=\begin{cases}t,&t\text{ even},\\t-t^{-1},&t\text{ odd},\end{cases}
\]
with $\omega(0)=0$, and it is attained when $|r-s|\leq1$, because $\psi(r,s)=t-(r-s)^2/t$. Dividing the runs equally within each pair is thus immediate. The question is how the totals themselves should be chosen.

\begin{theorem}\label{thm:nuisancecompleteclass}
Assume $n\geq2q$.
\begin{enumerate}
\item[(i)] If $n$ is even, every $D$-optimal design for~\eqref{eq:generalnuisanceJ} has $r_i=s_i$ for all $i$.
\item[(ii)] If $n$ is odd, every $D$-optimal design has $r_i=s_i$ for all pairs except one. On that pair $r_i$ and $s_i$ are positive and differ by one, so the design consists of an $(n-1)$-run design with $r_i=s_i$ for all $i$ and one further replicate of one of its runs.
\item[(iii)] For $n$ even or odd and any fixed nonnegative definite matrix $W$, there is a design of the form described in (i) or (ii) that minimises $\tr(WJ^{-1})$ among all designs with positive definite $J$.
\end{enumerate}
The design in (iii) may depend on $W$. In (ii), the $(n-1)$-run design need not be $D$-optimal for $n-1$ runs; an example is given in Section~S8 of the Supplementary Material.
\end{theorem}
\begin{proof}
Dividing the runs equally within a pair with a given total increases the information matrix. The remaining question is whether two pairs can have odd totals in an optimal design. For odd totals $t_i,t_j$ at two pairs with vectors $f_i,f_j$, write
\[
 J=A-\Gamma,\qquad J_\pm=A\pm\Theta,\qquad
 \Theta=f_if_i^\top-f_jf_j^\top,\quad \Gamma=t_i^{-1}f_if_i^\top+t_j^{-1}f_jf_j^\top,
\]
where $A$ contains the contributions of the other pairs and $t_if_if_i^\top+t_jf_jf_j^\top$. The matrices $J_\pm$ are obtained by changing the totals to $(t_i+1,t_j-1)$ or $(t_i-1,t_j+1)$ and dividing the runs equally in both pairs, and they are positive definite whenever $J$ is. Since $\Theta A^{-1}\Theta\preceq\Gamma$, the identity for the harmonic mean of two matrices gives
\[
 \left\{\tfrac12(J_+^{-1}+J_-^{-1})\right\}^{-1}
 =A-\Theta A^{-1}\Theta\succeq J.
\]
Thus the average of the two covariance matrices is not larger than the original one. At least one of the two changes does not increase a given value of $\tr(WJ^{-1})$, and by the strict concavity of $\log\det$ one of them strictly increases the determinant. Repeating the step leaves no odd total when $n$ is even and exactly one when $n$ is odd. A pair with a single run contributes no information, and moving that run to a pair with both orders present increases the determinant. The details, and the extension to losses that are concave and nondecreasing functions of the covariance matrix, are given in Section~S8 of the Supplementary Material.
\end{proof}

For $n=2h$ it follows that
\begin{equation}\label{eq:exactpairedreduction}
 \max_{\sum_i(r_i+s_i)=2h}\det J
 =2^q\max_{\substack{k_i\in\mathbb Z_{\geq0}\\\sum_i k_i=h}}
 \det\left(\sum_i k_i f_if_i^\top\right).
\end{equation}
Nothing is lost under the determinant criterion by considering only designs that contain both points of every pair equally often, although all designs are admitted as competitors, and the minimum of $\tr(WJ^{-1})$ is one half of the minimum of $\tr\{W(\sum_ik_if_if_i^\top)^{-1}\}$. A corresponding reduction for odd $n$ is given in Section~S8 of the Supplementary Material.

\begin{corollary}\label{cor:firstoddnuisance}
At $n=2q+1$, every $D$-optimal design for~\eqref{eq:generalnuisanceJ} is obtained by taking a $D$-optimal design with $2q$ runs and $r_i=s_i$ for all $i$, and replicating any one of its runs. Conversely, every design obtained in this way is $D$-optimal. The largest value of $\det J$ at $2q+1$ runs is exactly $4/3$ times the largest value at $2q$ runs.
\end{corollary}
\begin{proof}
By Theorem~\ref{thm:nuisancecompleteclass}, an optimal design consists of a $2q$-run design with $r_i=s_i$ for all $i$ and one further run on one of its pairs. Since positive definiteness requires at least $q$ pairs with both points present, the $2q$-run design has exactly $q$ distinct pairs, each with one run at each point, and their regression vectors are linearly independent. Let $\Phi_0$ be the $q\times q$ matrix with rows $f_1^\top,\ldots,f_q^\top$. The $2q$-run design has determinant $2^q(\det\Phi_0)^2$. The replicate changes one coefficient in~\eqref{eq:generalnuisanceJ} from $2$ to $8/3$, which multiplies the determinant by $4/3$ whichever run is replicated. Maximising over the choice of the $q$ independent vectors proves both directions.
\end{proof}

For the PWO factors, $\rho a=\rev a$ and $f=z$, and Theorems~\ref{thm:nuisancecompleteclass} and~\ref{thm:fixedsize} together give the following.

\begin{theorem}\label{thm:foldoveraug}
Let $m\geq3$ and let $D_0$ be a foldover OofA-OA with $n_0$ runs. Under model~\eqref{eq:nuisance}, for every order $a$ the design $D_0\cup\{a,\rev a\}$ maximises $\det\mathcal I_\beta$ and minimises $\tr(\Sigma\mathcal I_\beta^{-1})$ among all designs with $n_0+2$ runs.
\end{theorem}
\begin{proof}
By~\eqref{eq:exactpairedreduction}, $\det\mathcal I_\beta$ is maximised by a foldover design with $n_0/2+1$ pairs that maximises $\det G_D$, and by Theorem~\ref{thm:fixedsize} the design $D_0\cup\{a,\rev a\}$ is such a design. For the trace, Theorem~\ref{thm:nuisancecompleteclass}(iii) allows us to consider foldover designs only. For these $\mathcal I_\beta=nG_D/\sigma^2$ and $\tr(\Sigma G_D^{-1})=\tr(C_D^{-1})$, which is minimised by the same design according to Theorem~\ref{thm:fixedsize}.
\end{proof}

Theorem~\ref{thm:foldoveraug} compares designs for the estimation of $\beta$ after $g$ has been eliminated. Under the PWO model, where only the intercept is added to the PWO effects, the restriction to foldover designs in Theorem~\ref{thm:mainaugmentation}(ii) cannot be removed, as Section~\ref{sec:cost} shows.

Parts (i) and (ii) of Theorem~\ref{thm:nuisancecompleteclass} say that every $D$-optimal design has the stated form, whereas part (iii) only says that some optimal design has it. This cannot be improved. A design minimising $\tr(WJ^{-1})$ need not have $r_i=s_i$, and the reduction fails for the smallest eigenvalue; examples with three components are given in Section~S8 of the Supplementary Material. Related results are the approximate designs resistant to nuisance effects of \citet{rosaharman2016} and the optimal paired choice block designs of \citet{singhdaschai2019}, neither of which covers exact designs with an unrestricted number of runs in each pair.

The assumptions of Theorem~\ref{thm:nuisancecompleteclass} are met in other settings. In a two-level factorial experiment with factors coded $\pm1$, the map $x\mapsto-x$ pairs every run with its mirror image. The main effects and the interactions of an odd number of factors change sign under it, and the interactions of an even number of factors do not. Let $f(x)$ contain the main effects, with or without some interactions of an odd number of factors, and let the interactions of an even number of factors be part of an unrestricted term $g$. Theorem~\ref{thm:nuisancecompleteclass} then says that a $D$-optimal design consists of complete mirror-image pairs when $n$ is even, and of complete pairs and one replicate when $n$ is odd. This is the foldover principle of \citet{boxhunter1961} in the form of an optimality statement.

\subsection{The cost of the foldover restriction}\label{sec:cost}

Under the PWO model, adding a foldover pair to an OofA-OA is never the best way of adding two runs. For two orders $a$ and $b$ put
\[
 c(a,b)=x(a)^\top M_F^{-1}x(b)=1+\kappa(a,b),
\]
and define
\[
 g_m=\begin{cases}6,&m\text{ even},\\12,&m\text{ odd},\end{cases}
 \qquad
 \delta_m=\frac{\operatorname{dist}((m+1)p,g_m\mathbb Z)}{m+1},
\]
where $\operatorname{dist}(t,g_m\mathbb Z)$ is the distance from the integer $t$ to the nearest multiple of $g_m$.

\begin{proposition}\label{prop:tworun}
Let $D_0$ be an OofA-OA with $n_0$ runs on $m\geq3$ components, and let $D=D_0\cup\{a,b\}$, where $a$ and $b$ may coincide and may be in $D_0$. Put $n=n_0+2$ and $c=c(a,b)$.
\begin{enumerate}
\item[(i)] The design $D$ has
\begin{align}
 \frac{\det M_D}{\det M_F}
 &=\frac{n_0^{p-2}\{(n_0+p)^2-c^2\}}{n^p},
 \label{eq:tworundet}\\
 I(D)&=n\left\{\frac{p-2}{n_0}
       +\frac{2(n_0+p)}{(n_0+p)^2-c^2}\right\},
 \label{eq:tworunI}\\
 \lambda_{\min}(\mathcal C_D)&=\frac{n_0}{n},\qquad
 \tr\{(\mathcal C_D-I_p)^2\}
 =\frac{2p(p-2)+2c^2}{n^2}.
 \label{eq:tworunother}
\end{align}
Consequently the additions of two runs that minimise $|c(a,b)|$ maximise the determinant and minimise $I(D)$ and the squared distance at the same time, and all additions of two runs have the same minimum relative efficiency.
\item[(ii)] Replacing $\rev a$ in the addition $\{a,\rev a\}$ by an order obtained from $a$ by interchanging two adjacent components increases the determinant and decreases $I(D)$ and the squared distance. In particular, a foldover OofA-OA with an added foldover pair is not $D$-optimal among all designs of its run size when repeated runs are allowed.
\item[(iii)] For every two orders, $(m+1)c(a,b)\in(m+1)p+g_m\mathbb Z$ and $|c(a,b)|\geq\delta_m$. Hence every addition of two runs satisfies
\[
 \frac{\det M_D}{\det M_F}
 \leq\frac{n_0^{p-2}\{(n_0+p)^2-\delta_m^2\}}{(n_0+2)^p}.
\]
For $3\leq m\leq30$, equality holds, and $I(D)$ and the squared distance are smallest among additions of two runs, when $a$ is the identity order and $b$ is the order given in Table~S5 of the Supplementary Material.
\end{enumerate}
\end{proposition}

The proof is given in Section~S10 of the Supplementary Material. Part (i) follows from the eigenvalues $p\pm|c(a,b)|$ of $M_F^{-1/2}\{x(a)x(a)^\top+x(b)x(b)^\top\}M_F^{-1/2}$, and part (iii) from~\eqref{eq:kernel}. For two runs removed from an OofA-OA, \citet[Section~4.2]{tsai2022} gives the corresponding expression, with $(n_0-p)^2-c^2$ in place of $(n_0+p)^2-c^2$ and $n=n_0-2$, and removes the two runs with the smallest $|c(a,b)|$. The proposition compares additions to a given OofA-OA. It does not identify the optimal design among all designs with $n_0+2$ runs.
For the 24-run array generated by~\eqref{eq:base5}, the orders $12345$, $54321$ and $12543$ are not in the array. Adding the foldover pair $\{12345,54321\}$ gives $c=-9$, whereas adding $\{12345,12543\}$ gives $c=1$. Both designs have 26 distinct runs, and the values of $\det(X_D^\top X_D)$ are
\[
 66\,331\,474\,919\,424\quad\hbox{and}\quad
 70\,970\,039\,599\,104,
\]
with ratio $153/143$. The first design is optimal among foldover designs by Theorem~\ref{thm:mainaugmentation}. The second is optimal among additions of two runs to this array by part (iii), since $|c|=1=\delta_5$, but we do not claim that it is optimal among all 26-run designs. The pairs attaining the bound in (iii) may repeat runs of the array. When the added orders have to be new, the bound remains valid, and it is attained whenever two orders outside the array have the required value of $c(a,b)$, as in this example. For $m>30$ we do not know whether the bound is attained.

\subsection{Five components and 20 runs}\label{sec:twenty}

With five components $q=10$, so a foldover design in which all PWO effects are estimable needs at least 20 runs, and~\eqref{eq:fold} gives $\det M_D=(\det Z_H)^2/10^{10}$. No OofA-OA has 20 runs, and Theorem~\ref{thm:fixedsize} does not apply because $h=10$ is even. A repeated pair would make $Z_H$ singular when $h=q$, so the following result covers all foldover designs with 20 runs.

\begin{theorem}\label{thm:global5}
Among all foldover designs with five components and 20 runs,
\[
 \max_D\det M_D=\frac{2^{16}}{5^{10}},\qquad
 \max_D\Deff(D) \approx 0.897877 .
\]
Up to relabelling of the components there are exactly five nonequivalent designs attaining the maximum, and the eigenvalues of their matrices $C_D$ are of the two types shown in Table~\ref{tab:spectral5}.
\end{theorem}

\begin{table}[tb]
\centering
\caption{The five $D$-optimal foldover designs with five components and 20 runs. The multiplicity of each eigenvalue is given in parentheses.}
\label{tab:spectral5}
\begin{tabular}{@{}lcc@{}}
\toprule
 & Designs 1--3 & Designs 4--5\\
\midrule
Eigenvalues of $C_D$ & $2/5\ (3),\ 6/5\ (6),\ 8/5\ (1)$ & $3/5\ (6),\ 8/5\ (4)$\\
$I(D)$ & $14.125$ & $13.5$\\
$I(D)-1=\tr(C_D^{-1})$ & $13.125$ & $12.5$\\
$\eta(D)$ & $2/5$ & $3/5$\\
$\theta(D)$ & $\sqrt{7/375}$ & $\sqrt{2/75}$\\
Each component four times in each position & No & Yes\\
\bottomrule
\end{tabular}
\end{table}

The proof, which is a complete search in exact integer arithmetic, is given in Section~S11 of the Supplementary Material, and the five designs are listed in Section~S4. Counting relabellings separately, there are 276 optimal designs. Design 5 has a simple description as a set of affine maps over the field with five elements. The five designs have the same determinant but not the same moment matrix. Designs 4 and 5 have the smaller $I(D)$ and the larger minimum relative efficiency, and every component appears equally often in every position, whereas designs 1 to 3 are closer to the full design in the sense of $\theta(D)$. An experimenter who cares about the average variance of prediction, about the least favourable direction or about balance of the positions will prefer designs 4 and 5. We do not claim that they are optimal under these criteria among all 20-run foldover designs.

Theorem~\ref{thm:global5} is a statement about foldover designs. A search over all 120 orders finds a 20-run design that contains no foldover pair and has $\Deff\approx0.970089$, the same value as the 20-run design found by \citet{voelkel2019} by computer search (Section~S4 of the Supplementary Material); \citet[Table~6]{tsai2022} reports $0.9701$ for this run size. The restriction to foldover designs therefore costs at least this difference in efficiency under the PWO model, and we do not know the optimal 20-run design. Under model~\eqref{eq:nuisance}, however, that design has $\mathcal I_\beta=0$, because it contains no foldover pair, whereas the five foldover designs allow $\beta$ to be estimated. Theorem~\ref{thm:nuisancecompleteclass} then describes the optimal designs completely. Since the largest value of $(\det Z_H)^2$ for ten foldover pairs is $2^{26}$, the five designs are the $D$-optimal 20-run designs under~\eqref{eq:nuisance}, with $\det J=2^{36}$, and by Corollary~\ref{cor:firstoddnuisance} every $D$-optimal 21-run design replicates one run of one of them, which gives $\det J=2^{38}/3$.

\subsection{A sensitivity comparison}\label{sec:sensitivity}

Model~\eqref{eq:nuisance} gives a reason for using foldover pairs, but it does not say when their protection is worth the loss of efficiency under the PWO model. We compare designs with the same run size when the PWO model is fitted and the true mean is $\beta_0+z(a)^\top\beta+\phi(a)^\top\gamma$. The $r$ regressors in $\phi(a)$ have mean zero and variance one and are uncorrelated with the PWO factors over the full design, and $\gamma$ is random with mean zero and covariance matrix $\tau^2I_r/r$, so that $\tau/\sigma$ compares the size of the omitted effects with the error standard deviation. The mean squared error of the least squares estimator $\widehat\beta$ is
\[
 R_D(\tau)=E\{(\widehat\beta-\beta)^\top\Sigma(\widehat\beta-\beta)\}
 =\sigma^2\mathcal V_D+\tau^2\mathcal B_D,
\]
where $\mathcal V_D$ comes from the variance and $\mathcal B_D$ from the squared bias; the formulas are given in Section~S12 of the Supplementary Material.

For five components we consider two families of omitted effects, one at a time. The first consists of the ten standardised indicators that two given components are next to each other in the order. They do not change when the order is reversed, so that their contribution $\mathcal B_D^+$ to the bias is zero for a foldover design. The second consists of the thirty products of such an indicator with the PWO factor of two other components. They change sign under reversal, and a foldover design gives no protection against them. Table~\ref{tab:structuredrisk} gives the coefficients for four designs, namely design 5 of Table~\ref{tab:spectral5} and the 20-run design without foldover pairs of Section~\ref{sec:twenty}, and the array~\eqref{eq:base5} with $\{12345,54321\}$ or with $\{12345,12543\}$ added.

\begin{table}[tb]
\centering
\caption{Five components: $D$-efficiency under the PWO model and the coefficients in $R_D(\tau)/\sigma^2=\mathcal V_D+\mathcal B_D(\tau/\sigma)^2$. The column $\mathcal B_D^+$ refers to the omitted effects that are unchanged by reversal and $\mathcal B_D^-$ to those that change sign.}
\label{tab:structuredrisk}
\begin{tabular}{@{}rlrrrr@{}}
\toprule
$n$ & Design & $\Deff$ & $\mathcal V_D$ & $\mathcal B_D^+$ & $\mathcal B_D^-$\\
\midrule
20 & Foldover design 5 & 0.897877 & 0.625000 & 0 & 1.430556\\
20 & Design without foldover pairs & 0.970089 & 0.532417 & 0.462268 & 0.344876\\
26 & Array with $\{12345,54321\}$ & 0.982492 & 0.397727 & 0 & 0.816230\\
26 & Array with $\{12345,12543\}$ & 0.988548 & 0.392838 & 0.017893 & 0.753301\\
\bottomrule
\end{tabular}
\end{table}

For the first family, the foldover design has the smaller mean squared error exactly when $\tau/\sigma$ exceeds $0.447528$ at 20 runs and $0.522745$ at 26 runs. These thresholds belong to the designs and the family of effects considered here and are not general rules. For the second family the designs that are not foldover designs have the smaller variance and the smaller bias, hence the smaller mean squared error for every $\tau$. Foldover designs thus protect against omitted effects with one particular symmetry and not against a misspecified model in general.

\section{Concluding remarks}\label{sec:discussion}

Run sizes one above that of an OofA-OA are now settled, and so are foldover designs with one pair more than a foldover OofA-OA. Run sizes one below are not. The proofs go through the distance from the full design, and it is worth recording what this distance does not do. An OofA-OA with a run deleted is as close to the full design as its run size allows, and in the cases of Table~\ref{tab:deletioncomparison} it is not $D$-optimal. Away from the run sizes of Theorem~\ref{thm:mainaugmentation}, a moment matrix close to that of the full design and a large determinant are different aims.

We do not know for which numbers of components and run sizes an OofA-OA with one run deleted is $D$-optimal, nor which designs are optimal in the cases of Table~\ref{tab:deletioncomparison}. For two added runs we know the optimal design among all designs only for three components and for four components with 14 runs (Sections~S7 and~S10 of the Supplementary Material). For more than four components the bounds $B_m(h)$ need not be attained, because the conditions for equality on different sets of five labels may be incompatible, and the optimal designs at most run sizes remain unknown. We also do not know whether a foldover design with nine components and 74 runs attains the bounds of Theorem~\ref{thm:fixedsize}, or whether a design with ten components and 49 runs attains those of Corollary~\ref{cor:ordinaryfixed}.

Whether a foldover design should be used depends on which departures from the PWO model are feared. When effects that are unchanged by reversal are plausible and not small, Theorem~\ref{thm:nuisancecompleteclass} points to foldover designs, at a loss of efficiency under the PWO model that can be computed beforehand. For effects that change sign under reversal, the designs that are not foldover designs had the smaller mean squared error in the comparison of Section~\ref{sec:sensitivity}, but that comparison covers two run sizes and averages over families of effects, and for a single omitted effect a design chosen for it can do better than a design of Theorem~\ref{thm:mainaugmentation}(i) (Section~S12 of the Supplementary Material). When no effects outside the PWO model are expected, the designs of Theorem~\ref{thm:mainaugmentation}(i) are optimal at their run sizes.

\section*{Supplementary Material}

The online Supplementary Material contains the proofs that are not given in the text, the designs, the enumerations and the details of the sensitivity comparison.

\vspace{12pt}
\noindent Stelios D. Georgiou and Stella Stylianou\\
Department of Mathematical and Geospatial Sciences, School of Science, RMIT University, Melbourne, Australia\\
E-mail: stelios.georgiou@rmit.edu.au; stella.stylianou@rmit.edu.au

\vspace{8pt}
\noindent Haralambos Evangelaras\\
Department of Statistics and Insurance Science, University of Piraeus, Piraeus, Greece\\
E-mail: hevangel@unipi.gr
\end{document}


\begin{center}
{\large\bfseries Supplementary Material for ``Exact optimality and foldover designs in order-of-addition experiments''}\\[14pt]
Stelios D. Georgiou$^{1}$, Haralambos Evangelaras$^{2}$ and Stella Stylianou$^{1}$\\[8pt]
{\itshape $^{1}$RMIT University and $^{2}$University of Piraeus}
\end{center}

This supplement contains the proofs that are not given in the paper, the designs used in its examples and several enumerations. Numbers with the prefix S refer to this supplement, and all other numbers refer to the paper. Orders are written as strings of component labels in the order of addition. For a foldover design only the half-design $H$ is listed, and the design consists of the listed orders and their reverses. We call the listed order the representative of its foldover pair, and a multiset is a collection in which elements may be repeated.

We keep the notation of the paper: $q=m(m-1)/2$ is the number of PWO factors, $p=q+1$, and a foldover design $D=H\cup\rev H$ has $n=2h$ runs, where $h$ counts the foldover pairs with multiplicity. Its matrices are $G_D=Z_H^\top Z_H/h$ and $C_D=\Sigma^{-1/2}G_D\Sigma^{-1/2}$. A matrix $Z$ without subscript is the PWO matrix of the design under discussion. The criteria satisfy
\[
\begin{gathered}
 V(H)=h^2\Delta(D)=h^2q(q-1)\theta(D)^2,\\
 I(D)=1+\tr(C_D^{-1}),\qquad \Deff(D)=(\det C_D)^{1/p}.
\end{gathered}
\]
Efficiencies compare moment matrices, that is, information per run, so a larger design can have a smaller efficiency although it carries more information in total. The symbol $h_0$ refers to a given OofA-OA and $h$ to the final number of foldover pairs.

\section{Designs}\label{sec:designs}

\noindent\textit{Five components.} The 24-run foldover OofA-OA used in Section~3 consists of the twelve orders in~(3.1) of the paper and their reverses. Its half-design has $Z_H^\top Z_H=12\Sigma$, and on every set of four components each of the twelve foldover pairs occurs exactly once in $H$, in agreement with the case $\theta=0$ discussed after Proposition~1.

\noindent\textit{Six components.} The 48-run foldover OofA-OA used in Section~3 and in Section~\ref{sec:mod} has the half-design
\[
\begin{array}{llllll}
\mathtt{123456}&\mathtt{126543}&\mathtt{135624}&\mathtt{142365}&\mathtt{146352}&\mathtt{153426}\\
\mathtt{156423}&\mathtt{163425}&\mathtt{165234}&\mathtt{214563}&\mathtt{235164}&\mathtt{241365}\\
\mathtt{254163}&\mathtt{261354}&\mathtt{263415}&\mathtt{312546}&\mathtt{321645}&\mathtt{341625}\\
\mathtt{354216}&\mathtt{362154}&\mathtt{421536}&\mathtt{431526}&\mathtt{462315}&\mathtt{523146}\\
\end{array}
\]
and consists of these orders and their reverses. All 48 runs are distinct, $3Z_H^\top Z_H=24(I_{15}+B^\top B)$, so $G_D=\Sigma$, and on every set of four components each foldover pair occurs twice in $H$. The design was found by an exchange search over foldover pairs that maximises the determinant in~(2.2). The PWO columns of the displayed orders do not all sum to zero, but replacing some of the orders by their reverses gives a 24-run OofA-OA: reading the display row by row, reverse the orders in positions $2,3,5,6,8,11,12,15,16,21,23,24$. The resulting matrix $Z_*$ satisfies
\[
 \mathbf1^\top Z_*=0,\qquad 3Z_*^\top Z_*=24(I_{15}+B^\top B).
\]
Thus the 48-run design consists of a 24-run OofA-OA and the reverses of its runs. The enumeration of the OofA-OAs with 24 runs by \citet{schoenmee2023} is complete, so this array is equivalent to one of their $22\,651$ nonequivalent arrays with six components, but we have not determined to which. By Corollary~2, 48 is the smallest possible run size of a foldover OofA-OA with six components.

\noindent\textit{Seven components.} The 24-run OofA-OA with seven components used in Table~2 and in Section~5 is the array of \citet[Example~5]{tsai2022}, which is also the array with seven components and 24 runs in the \emph{OofAExp} catalogue, version~0.1.0 \citep{tsai2026package}. Its runs are
\[
\begin{array}{llllll}
\mathtt{1245376}&\mathtt{1643572}&\mathtt{1726534}&\mathtt{1735624}&\mathtt{2154673}&\mathtt{2356471}\\
\mathtt{2731645}&\mathtt{2746135}&\mathtt{3265174}&\mathtt{3461572}&\mathtt{3742156}&\mathtt{3751246}\\
\mathtt{4316275}&\mathtt{4512673}&\mathtt{4753261}&\mathtt{4762351}&\mathtt{5421376}&\mathtt{5623174}\\
\mathtt{5713462}&\mathtt{5764312}&\mathtt{6134275}&\mathtt{6532471}&\mathtt{6715423}&\mathtt{6724513}\\
\end{array}
\]
Its 24 runs are distinct, $\mathbf1^\top Z=0$ and $Z^\top Z=24\Sigma$. Together with the reverses of its runs it forms the 48-run foldover OofA-OA of Section~5, and adding any order that is not in the array gives the optimal 25-run design of Table~2.

\noindent\textit{Eight components.} The 48-run OofA-OA of Table~2 is the array with eight components and 48 runs in the same version of the catalogue, which attributes it to \citet{tsai2022}. It is generated from the six orders $\mathtt{02514376}$, $\mathtt{04135267}$, $\mathtt{04172653}$, $\mathtt{06352417}$, $\mathtt{06437512}$, $\mathtt{06712453}$ by adding $b=0,\ldots,7$ modulo eight to every entry and then adding one to every label. The resulting runs are
\[
\begin{array}{llllll}
\mathtt{13625487}&\mathtt{24736518}&\mathtt{35847621}&\mathtt{46158732}&\mathtt{57261843}&\mathtt{68372154}\\
\mathtt{71483265}&\mathtt{82514376}&\mathtt{15246378}&\mathtt{26357481}&\mathtt{37468512}&\mathtt{48571623}\\
\mathtt{51682734}&\mathtt{62713845}&\mathtt{73824156}&\mathtt{84135267}&\mathtt{15283764}&\mathtt{26314875}\\
\mathtt{37425186}&\mathtt{48536217}&\mathtt{51647328}&\mathtt{62758431}&\mathtt{73861542}&\mathtt{84172653}\\
\mathtt{17463528}&\mathtt{28574631}&\mathtt{31685742}&\mathtt{42716853}&\mathtt{53827164}&\mathtt{64138275}\\
\mathtt{75241386}&\mathtt{86352417}&\mathtt{17548623}&\mathtt{28651734}&\mathtt{31762845}&\mathtt{42873156}\\
\mathtt{53184267}&\mathtt{64215378}&\mathtt{75326481}&\mathtt{86437512}&\mathtt{17823564}&\mathtt{28134675}\\
\mathtt{31245786}&\mathtt{42356817}&\mathtt{53467128}&\mathtt{64578231}&\mathtt{75681342}&\mathtt{86712453}\\
\end{array}
\]

\noindent\textit{Nine components.} The 48-run OofA-OA with nine components given by \citet{tsai2025} is
\[
\begin{array}{lll}
\mathtt{136295487}&\mathtt{152496378}&\mathtt{152837964}\\
\mathtt{179548623}&\mathtt{247396518}&\mathtt{263148975}\\
\mathtt{263597481}&\mathtt{281346759}&\mathtt{285746319}\\
\mathtt{289651734}&\mathtt{312457869}&\mathtt{316857429}\\
\mathtt{319762845}&\mathtt{358497621}&\mathtt{374251986}\\
\mathtt{374698512}&\mathtt{429873156}&\mathtt{461598732}\\
\mathtt{485362917}&\mathtt{485791623}&\mathtt{516473928}\\
\mathtt{516892734}&\mathtt{539184267}&\mathtt{572691843}\\
\mathtt{627193845}&\mathtt{627584931}&\mathtt{641382759}\\
\mathtt{645782319}&\mathtt{649215378}&\mathtt{683792154}\\
\mathtt{714893265}&\mathtt{738294156}&\mathtt{738615942}\\
\mathtt{752413869}&\mathtt{756813429}&\mathtt{759326481}\\
\mathtt{825194376}&\mathtt{841395267}&\mathtt{841726953}\\
\mathtt{869437512}&\mathtt{917463528}&\mathtt{917823564}\\
\mathtt{942356817}&\mathtt{942716853}&\mathtt{953467128}\\
\mathtt{953827164}&\mathtt{986352417}&\mathtt{986712453}\\
\end{array}
\]
Its runs are distinct, $\mathbf1^\top Z=0$ and $Z^\top Z=48\Sigma$. The order $123456789$ is not in the array, and adding it gives the optimal 49-run design of Table~2.

\noindent\textit{Ten components.} Identify the labels $1,\ldots,10$ with the unordered pairs
\[
\{1,2\},\ \{1,3\},\ \{1,4\},\ \{1,5\},\ \{2,3\},\ \{2,4\},\ \{2,5\},\ \{3,4\},\ \{3,5\},\ \{4,5\}
\]
of five auxiliary symbols, in this order, so that a permutation $g$ of the five symbols relabels the component $\{i,j\}$ as $\{g(i),g(j)\}$. The two orders $c_1=(1,2,3,7,8,9,6,4,10,5)$ and $c_2=(1,2,3,6,8,7,9,10,5,4)$ have trivial stabilisers, so each generates an orbit of 120 distinct runs. The two orbits are disjoint, and their union $D_0$ satisfies $\mathbf1^\top Z_{D_0}=0$ and $Z_{D_0}^\top Z_{D_0}=240\Sigma$, as checked in exact integer arithmetic. The identity order is not in $D_0$.

\section{The construction in Proposition~2}\label{sec:nuisanceproof}

Reverse $(1,\ldots,m)$ by $q$ interchanges of adjacent components, in which each pair of components is interchanged exactly once: move $m$ to the first position, then $m-1$ to the second, and so on. Write $a^{(r)}$ for the order after step $r$ and $z^{(r)}=z(a^{(r)})$, so that superscripts count steps and the subscripts $jk$ still refer to pairs of components. Each step changes one coordinate by $-2$. When the coordinates are arranged in the order of the interchanges, the differences $z^{(r)}-z^{(0)}$, $r=1,\ldots,q$, form a triangular matrix with diagonal entries $-2$. They are linearly independent, so $z^{(0)},\ldots,z^{(q)}$ are affinely independent and span $\R^q$. Since $z^{(q)}=-z^{(0)}$, their only linear relation is a multiple of $z^{(0)}+z^{(q)}=0$. A relation among the first $q$ vectors would have coefficient zero on $z^{(q)}$, which forces every coefficient to vanish. Thus $z^{(0)},\ldots,z^{(q-1)}$ are linearly independent and contain no foldover pair. These orders and their reverses give $2q$ distinct runs, and the PWO effects are estimable under model~(6.1) with the smallest possible number of runs.

Two limits of the protection given by foldover pairs should be noted. With $h=q$ pairs the regression of the half-differences on the PWO factors is saturated, so further pairs or replicates are needed to estimate the error variance. The elimination of $g$ protects the estimators of the PWO effects against every reversal-invariant term but not against omitted effects that change sign under reversal. Moreover, the criterion $I(D)$ of the paper concerns the PWO part of the mean, since an unrestricted $g(a)$ at a pair that has not been run cannot be predicted from the differences within pairs.

\section{Weaker bounds and enumeration for four components}\label{sec:weaker}

Lemma~3 of the paper uses two contributions to the squared norm in~(4.5): a contribution from the parity of the entries for pairs of triples sharing at most one label, available when $h$ is odd, and a contribution from the four-element sets, available when $h$ is not divisible by three. When only one of them is available the same argument gives the following weaker bounds.

\begin{lemma}\label{lem:weaker}
Let $H$ be any multiset of $h$ orders on $m\geq3$ components and let $V(H)$ be as in~(4.5).
\begin{enumerate}
\item[(i)] If $h$ is even and not divisible by three, then $V(H)\geq4m(m-1)(m-2)/(m+1)$.
\item[(ii)] If $h$ is odd and $3\mid h$, then $V(H)\geq q(q-1)-4m(m-1)(m-2)/(m+1)$.
\end{enumerate}
For the foldover design $D=H\cup\rev H$, these bounds give $h\theta(D)\geq c$, with $c=\{8(m-2)/((m+1)(q-1))\}^{1/2}$ and $c=\{1-8(m-2)/((m+1)(q-1))\}^{1/2}$, respectively.
\end{lemma}
\begin{proof}
In the proof of Lemma~3, under (i) only the contribution of the four-element sets is available, giving $\|U^\top U-hR\|_F^2\geq\tfrac{32}3\binom{m+1}4$, and under (ii) only the parity contribution, giving $|\mathcal T|(|\mathcal T|-1)-12\binom{m+1}4$; multiplying by $9/(m+1)^2$ gives the stated values, since $9\cdot\tfrac{32}{3}\binom{m+1}4/(m+1)^2=4m(m-1)(m-2)/(m+1)$.
\end{proof}

Neither bound is attained in the cases with $m\geq4$ enumerated below. Theorem~4 gives stronger bounds for these values of $h$ and a positive bound when $h\equiv6\pmod{12}$.

\begin{table}[tb]
\centering
\caption{Four components, all multisets of $h$ foldover pairs. The column headed bound is $\{((h-1)/h)^{5}(h+5)/h\}^{1/7}$; it is a proved bound for every odd $h>1$, namely $h=3,5,7,9,11$ (in bold), and is shown for reference otherwise.}
\label{tab:m4}
\begin{tabular}{rrrrrrrr}
\toprule
$h$ & multisets & $\max\det(Z_H^\top Z_H)$ & $\Deff$ & optima & bound & $\min V$ & $h\theta_{\min}$\\
\midrule
1 & 12 & 0 & 0 & -- & -- & 30.00 & 1.000\\
2 & 78 & 0 & 0 & -- & 0.7290 & 48.00 & 1.265\\
3 & 364 & 0 & 0 & -- & \textbf{0.8611} & 62.64 & 1.445\\
4 & 1365 & 0 & 0 & -- & 0.9143 & 65.28 & 1.475\\
5 & 4368 & 0 & 0 & -- & \textbf{0.9414} & 58.80 & 1.400\\
6 & 12376 & 4096 & 0.9088 & 2 & 0.9573 & 43.20 & 1.200\\
7 & 31824 & 10240 & 0.9077 & 12 & \textbf{0.9674} & 58.80 & 1.400\\
8 & 75582 & 24576 & 0.9174 & 30 & 0.9743 & 65.28 & 1.475\\
9 & 167960 & 57344 & 0.9360 & 20 & \textbf{0.9792} & 62.64 & 1.445\\
10 & 352716 & 128000 & 0.9591 & 6 & 0.9828 & 48.00 & 1.265\\
11 & 705432 & 256000 & 0.9759 & 12 & \textbf{0.9856} & 30.00 & 1.000\\
12 & 1352078 & 512000 & 1.0000 & 1 & 0.9877 & 0.00 & 0.000\\
\bottomrule
\end{tabular}
\end{table}

For four components there are twelve foldover pairs, and a foldover design with $h$ pairs, repeated pairs allowed, is a multiset of size $h$ from these twelve. Table~\ref{tab:m4} lists, for every $h\leq12$, the number of multisets, the maximum of $\det(Z_H^\top Z_H)$ and the corresponding efficiency, the number of multisets attaining the maximum, and the minimum of $V(H)$ together with $h\theta_{\min}$. For $h\geq6$, all maxima are attained only by designs with distinct runs. The small decrease in efficiency from $h=6$ to $h=7$ is due to the division by the number of runs and does not mean that information is lost when a pair is added. For $h<6$ every design is singular, and the column of optima is left blank. The bound is the right-hand side of the first inequality in~(4.9) of the paper. It is proved for every odd $h>1$ by Theorem~6, and at $h=12$ it is exceeded by the full design, which shows that the condition that $h$ is odd cannot be dropped. The minimum of $V(H)$ equals $q(q-1)=30$ exactly at $h=1$ and $h=11$, where the bound of Lemma~3 is attained. For the other values of $h$ covered by Lemma~\ref{lem:weaker}, $h\theta_{\min}$ exceeds the lower bounds $0.8$ and $0.6$ of that lemma by a wide margin. The optimal designs at $h=7$ and at $h=11$ each form one equivalence class of twelve designs. The two optimal designs at $h=6$ consist of the twelve even permutations and of the twelve odd permutations of the four components. Every optimal design at $h=7$ consists of one of these two designs and one further pair, and the optimal designs at $h=11$ are the full design with one pair deleted.

For five components the smallest values of $V(H)$ over all multisets with $h=2,3,4$ pairs are $160$, $210$ and $248$, that is $h\theta_{\min}\approx1.333$, $1.528$ and $1.660$, against the constants $0.667$ and $0.745$ of Lemma~\ref{lem:weaker}. The bounds of Theorem~4 are $B_5(2)=130$, $B_5(3)=158$ and $B_5(4)=136$, and $B_5(10)=130$ gives the bound $0.9626$ on the efficiency in Corollary~\ref{cor:allresiduecriteria} at 20 runs, against the optimal value $0.8979$ of Theorem~9. For more than four components these bounds are therefore valid but not sharp for small $h$. With $h<10$ the PWO effects are not estimable in a foldover design with five components, and $h=10$ is covered by Theorem~9.

\section{Enumeration and classification for Theorem~9}\label{sec:enum}

The search selects ten of the sixty foldover pairs on five components, each represented by its lexicographically smaller order, with the order $12345$ fixed by relabelling. It bounds the squared determinant of every completion by the following form of Hadamard's inequality in terms of Schur complements \citep{hornjohnson2013}.

\begin{lemma}\label{lem:branch}
Let $x_1,\ldots,x_r\in\R^q$ be independent, with Gram matrix $P$ and $\delta=\det P>0$. For an available vector $x_i$, put $b_i=\det\begin{pmatrix}P&c_i\\c_i^\top&x_i^\top x_i\end{pmatrix}$ with $c_i=(x_1^\top x_i,\ldots,x_r^\top x_i)^\top$. If $k=q-r$ further rows are required, the squared determinant of any completion is at most $b_{(1)}\cdots b_{(k)}/\delta^{k-1}$, where $b_{(1)}\geq\cdots\geq b_{(k)}$ are the largest available bordered determinants.
\end{lemma}
A set of dependent rows is discarded, since a square matrix containing dependent rows is singular. The inequality is standard and is stated here because the search uses it.

For this enumeration, let $Z$ be the $60\times10$ PWO matrix of all available representatives and $\Gamma=ZZ^\top$ its raw row Gram matrix. Let $E$ denote the matrix of scaled residual inner products. After fixing the identity row, indexed by $0$, one has $\delta=10$ and $E_{ij}=10\Gamma_{ij}-\Gamma_{i0}\Gamma_{0j}$; in general $E_{ij}=\delta\,x_i^\top(I-\Pi)x_j$, where $\Pi$ projects onto the selected row span. The entries $E_{ij}$ are bordered Gram minors and integers, and $b_i=E_{ii}$. Appending row $j$ replaces $\delta$ by $E_{jj}$ and updates $E'_{uv}=(E_{jj}E_{uv}-E_{uj}E_{jv})/\delta$, the division being exact by Sylvester's identity \citep{bareiss1968}.

The candidates are taken in increasing lexicographic order, so a node with $N$ remaining candidates and $k$ rows still required accounts for $\binom Nk$ completions, either by discarding the branch or through its children, and the accumulated count verifies that all $\binom{59}9=12\,565\,671\,261$ selections are covered. A branch is discarded only when the bound of Lemma~\ref{lem:branch} is strictly below $T=8192^2$, so that every selection attaining the maximum is retained. The search visits $5\,787\,915$ nodes, where a node is a recursive visit with at least two and at most nine selected rows; the common root with the identity row alone and the 327 complete ten-row leaves that survive pruning are not counted, a convention stated so that independent implementations can be compared.

Exactly 46 retained selections have $|\det Z_H|=8192$, each verified by integer elimination, and none has a larger determinant. Each of the 46 retained determinants was separately evaluated by scalar fraction-free elimination, and the five orbits under relabelling of the components were checked against the retained selections. They are listed in Table~\ref{tab:sel46}.

\begin{table}[tb]
\centering
\caption{The five nonequivalent $D$-optimal foldover designs with five components and 20 runs. Design $i$ consists of the ten orders in $H_i$ and their reverses; $c_i$ is the number of retained selections equivalent to it, $a_i$ the order of its stabiliser under relabelling of the components, and $120/a_i$ the number of designs equivalent to it. All five have $\Deff\approx0.8979$, and their eigenvalues are given in Table~4 of the paper.}
\label{tab:optimal5}
\begin{tabular}{clrrr}
\toprule
Design & $H_i$ & $c_i$ & $a_i$ & $120/a_i$\\
\midrule
1 & \(\begin{array}{lllll}\mathtt{12345} & \mathtt{12543} & \mathtt{13524} & \mathtt{14325} & \mathtt{15342}\\\mathtt{21453} & \mathtt{23154} & \mathtt{25134} & \mathtt{32415} & \mathtt{42135}\end{array}\) & 20 & 1 & 120\\[3pt]
2 & \(\begin{array}{lllll}\mathtt{12345} & \mathtt{12543} & \mathtt{13524} & \mathtt{14532} & \mathtt{21435}\\\mathtt{21534} & \mathtt{24513} & \mathtt{32415} & \mathtt{32514} & \mathtt{41325}\end{array}\) & 10 & 2 & 60\\[3pt]
3 & \(\begin{array}{lllll}\mathtt{12345} & \mathtt{12543} & \mathtt{13542} & \mathtt{14523} & \mathtt{21435}\\\mathtt{21534} & \mathtt{23514} & \mathtt{31425} & \mathtt{31524} & \mathtt{42315}\end{array}\) & 10 & 2 & 60\\[3pt]
4 & \(\begin{array}{lllll}\mathtt{12345} & \mathtt{13254} & \mathtt{14352} & \mathtt{15423} & \mathtt{21453}\\\mathtt{23514} & \mathtt{24135} & \mathtt{31524} & \mathtt{34215} & \mathtt{43125}\end{array}\) & 5 & 4 & 30\\[3pt]
5 & \(\begin{array}{lllll}\mathtt{12345} & \mathtt{13254} & \mathtt{14532} & \mathtt{15423} & \mathtt{21435}\\\mathtt{24153} & \mathtt{25314} & \mathtt{31524} & \mathtt{34215} & \mathtt{43125}\end{array}\) & 1 & 20 & 6\\[3pt]
\bottomrule
\end{tabular}
\end{table}

For the classification, every permutation in $S_5$ is applied to each retained selection, each order is replaced by the smaller order of its foldover pair, and the sorted ten orders are compared. The smallest image is taken as the canonical form, and the five distinct canonical forms are the designs of Table~\ref{tab:optimal5}, with stabilisers of orders $1,2,2,4,20$ and orbits of sizes $120,60,60,30,6$. The eigenvalues, the values of $I(D)$ and $\eta(D)$ and the position counts in Table~4 of the paper are computed directly from the five designs. Designs 1--3 have $I(D)=14.125$, that is $\tr(C_D^{-1})=13.125$, with $\tr(C_D^2)=11.68$ and $\theta\approx0.1366$. Designs 4 and 5 have $I(D)=13.5$, that is $\tr(C_D^{-1})=12.5$, with $\tr(C_D^2)=12.40$ and $\theta\approx0.1633$, and in these two designs every component is in every position exactly four times. Their smaller $I(D)$ does not show that they are $I$-optimal among all 20-run foldover designs.

A design equivalent to design 5 is the set of orders $(a\pi_1+b,\ldots,a\pi_5+b)$, $\pi=(0,1,4,2,3)$, $a\in\mathbb F_5\setminus\{0\}$, $b\in\mathbb F_5$, with arithmetic in the field $\mathbb F_5$ and labels converted to $1,\ldots,5$. It has twenty distinct runs and is a foldover design, because $\pi_j+\pi_{6-j}=3$ in $\mathbb F_5$ for every $j$, so that the reverse of the order indexed by $(a,b)$ is the order indexed by $(-a,3a+b)$. Its half-design has $|\det Z_H|=8192$.

Theorem~9 concerns foldover designs only. For comparison, an exchange search over all 120 orders was started from 5\,000 random sets of 20 runs, each drawn by a partial Fisher--Yates shuffle driven by the congruential generator of Section~\ref{sec:ranking} with seed 5. At each step it replaces one run by the exchange that increases the determinant most, until no exchange increases it. The search found the 20-run design
\[
\begin{array}{llllllllll}
\mathtt{12345}&\mathtt{12435}&\mathtt{13542}&\mathtt{15432}&\mathtt{25134}&\mathtt{25413}&\mathtt{31254}&\mathtt{32145}&\mathtt{34125}&\mathtt{34521}\\
\mathtt{35241}&\mathtt{41532}&\mathtt{42153}&\mathtt{42315}&\mathtt{43251}&\mathtt{45312}&\mathtt{51423}&\mathtt{52314}&\mathtt{52431}&\mathtt{53412}
\end{array}
\]
which contains no foldover pair and has $\det(X_D^\top X_D)=2^{24}\cdot3^3\cdot5\cdot7^2\cdot29$ and $\Deff\approx0.970089$, compared with $0.8979$ for the optimal foldover designs. In the recorded run this determinant was reached by 78 of the $5\,000$ restarts and exceeded by none, and the design shown is one of those obtained. The number of restarts reaching it depends on how ties are resolved in floating-point arithmetic. \citet{voelkel2019} reports two 20-run designs with five components found by computer search, with efficiencies $0.90$ under a chi-squared criterion and $0.97$ under the determinant criterion, and specifies the second by its row indices in the lexicographic ordering of the 120 orders. Recomputed from those indices, it has the same determinant as the design above. No optimality is claimed for either design. This comparison of efficiencies is under the PWO model. Under model~(6.1) of the paper, the information matrix of this design for the PWO effects is zero, whereas for each of the five optimal foldover designs it has full rank.

\FloatBarrier
\begingroup
\setlength{\LTcapwidth}{\linewidth}
\begin{longtable}{rl}
\caption{The 46 selections of ten orders containing $12345$ with $|\det Z_H|=8192$, one order from each foldover pair, in lexicographic order.}\label{tab:sel46}\\
\toprule
& selection\\
\midrule
\endfirsthead
\toprule
& selection\\
\midrule
\endhead
1&$\mathtt{12345}\ \mathtt{12543}\ \mathtt{13524}\ \mathtt{14325}\ \mathtt{15342}\ \mathtt{21453}\ \mathtt{23154}\ \mathtt{25134}\ \mathtt{32415}\ \mathtt{42135}$\\
2&$\mathtt{12345}\ \mathtt{12543}\ \mathtt{13524}\ \mathtt{14532}\ \mathtt{21435}\ \mathtt{21534}\ \mathtt{24513}\ \mathtt{32415}\ \mathtt{32514}\ \mathtt{41325}$\\
3&$\mathtt{12345}\ \mathtt{12543}\ \mathtt{13524}\ \mathtt{14532}\ \mathtt{21534}\ \mathtt{24315}\ \mathtt{24513}\ \mathtt{32514}\ \mathtt{34125}\ \mathtt{42135}$\\
4&$\mathtt{12345}\ \mathtt{12543}\ \mathtt{13542}\ \mathtt{14523}\ \mathtt{15324}\ \mathtt{21435}\ \mathtt{23154}\ \mathtt{25134}\ \mathtt{31425}\ \mathtt{35124}$\\
5&$\mathtt{12345}\ \mathtt{12543}\ \mathtt{13542}\ \mathtt{14523}\ \mathtt{21435}\ \mathtt{21534}\ \mathtt{23514}\ \mathtt{31425}\ \mathtt{31524}\ \mathtt{42315}$\\
6&$\mathtt{12345}\ \mathtt{12543}\ \mathtt{13542}\ \mathtt{14532}\ \mathtt{21435}\ \mathtt{25134}\ \mathtt{31425}\ \mathtt{32415}\ \mathtt{35124}\ \mathtt{41325}$\\
7&$\mathtt{12345}\ \mathtt{12543}\ \mathtt{13542}\ \mathtt{14532}\ \mathtt{21534}\ \mathtt{24135}\ \mathtt{31524}\ \mathtt{32514}\ \mathtt{34125}\ \mathtt{42315}$\\
8&$\mathtt{12345}\ \mathtt{12543}\ \mathtt{14325}\ \mathtt{15342}\ \mathtt{21354}\ \mathtt{21453}\ \mathtt{25314}\ \mathtt{31524}\ \mathtt{32415}\ \mathtt{42315}$\\
9&$\mathtt{12345}\ \mathtt{12543}\ \mathtt{14352}\ \mathtt{15324}\ \mathtt{21354}\ \mathtt{21453}\ \mathtt{24315}\ \mathtt{31425}\ \mathtt{32514}\ \mathtt{41325}$\\
10&$\mathtt{12345}\ \mathtt{12543}\ \mathtt{14352}\ \mathtt{15324}\ \mathtt{21354}\ \mathtt{24315}\ \mathtt{24513}\ \mathtt{32145}\ \mathtt{35124}\ \mathtt{41325}$\\
11&$\mathtt{12345}\ \mathtt{12543}\ \mathtt{14352}\ \mathtt{15342}\ \mathtt{21354}\ \mathtt{24153}\ \mathtt{31524}\ \mathtt{32145}\ \mathtt{41325}\ \mathtt{42315}$\\
12&$\mathtt{12345}\ \mathtt{12543}\ \mathtt{14352}\ \mathtt{15342}\ \mathtt{21453}\ \mathtt{23154}\ \mathtt{31425}\ \mathtt{32415}\ \mathtt{32514}\ \mathtt{42135}$\\
13&$\mathtt{12345}\ \mathtt{13254}\ \mathtt{14253}\ \mathtt{14352}\ \mathtt{15243}\ \mathtt{23514}\ \mathtt{24513}\ \mathtt{32145}\ \mathtt{42135}\ \mathtt{43125}$\\
14&$\mathtt{12345}\ \mathtt{13254}\ \mathtt{14253}\ \mathtt{15342}\ \mathtt{23154}\ \mathtt{24153}\ \mathtt{25134}\ \mathtt{32145}\ \mathtt{34125}\ \mathtt{42135}$\\
15&$\mathtt{12345}\ \mathtt{13254}\ \mathtt{14253}\ \mathtt{21534}\ \mathtt{23514}\ \mathtt{24315}\ \mathtt{24513}\ \mathtt{32145}\ \mathtt{34125}\ \mathtt{42135}$\\
16&$\mathtt{12345}\ \mathtt{13254}\ \mathtt{14352}\ \mathtt{15243}\ \mathtt{15342}\ \mathtt{21453}\ \mathtt{23154}\ \mathtt{32415}\ \mathtt{42135}\ \mathtt{43125}$\\
17&$\mathtt{12345}\ \mathtt{13254}\ \mathtt{14352}\ \mathtt{15243}\ \mathtt{23154}\ \mathtt{24135}\ \mathtt{25143}\ \mathtt{32145}\ \mathtt{35124}\ \mathtt{43125}$\\
18&$\mathtt{12345}\ \mathtt{13254}\ \mathtt{14352}\ \mathtt{15423}\ \mathtt{21453}\ \mathtt{23514}\ \mathtt{24135}\ \mathtt{31524}\ \mathtt{34215}\ \mathtt{43125}$\\
19&$\mathtt{12345}\ \mathtt{13254}\ \mathtt{14532}\ \mathtt{15243}\ \mathtt{21435}\ \mathtt{24513}\ \mathtt{25314}\ \mathtt{32415}\ \mathtt{35124}\ \mathtt{43125}$\\
20&$\mathtt{12345}\ \mathtt{13254}\ \mathtt{14532}\ \mathtt{15423}\ \mathtt{21435}\ \mathtt{24153}\ \mathtt{25314}\ \mathtt{31524}\ \mathtt{34215}\ \mathtt{43125}$\\
21&$\mathtt{12345}\ \mathtt{13254}\ \mathtt{14532}\ \mathtt{15423}\ \mathtt{21534}\ \mathtt{24315}\ \mathtt{25143}\ \mathtt{31425}\ \mathtt{35214}\ \mathtt{42135}$\\
22&$\mathtt{12345}\ \mathtt{13254}\ \mathtt{15342}\ \mathtt{21453}\ \mathtt{23154}\ \mathtt{25134}\ \mathtt{31425}\ \mathtt{32415}\ \mathtt{35214}\ \mathtt{42135}$\\
23&$\mathtt{12345}\ \mathtt{13524}\ \mathtt{14325}\ \mathtt{15423}\ \mathtt{21453}\ \mathtt{23514}\ \mathtt{24315}\ \mathtt{25413}\ \mathtt{31254}\ \mathtt{34215}$\\
24&$\mathtt{12345}\ \mathtt{13524}\ \mathtt{14325}\ \mathtt{15423}\ \mathtt{21534}\ \mathtt{23514}\ \mathtt{24315}\ \mathtt{25413}\ \mathtt{34125}\ \mathtt{42135}$\\
25&$\mathtt{12345}\ \mathtt{13524}\ \mathtt{14523}\ \mathtt{15342}\ \mathtt{23154}\ \mathtt{24153}\ \mathtt{25143}\ \mathtt{32145}\ \mathtt{42135}\ \mathtt{43125}$\\
26&$\mathtt{12345}\ \mathtt{13524}\ \mathtt{14532}\ \mathtt{15423}\ \mathtt{21435}\ \mathtt{24153}\ \mathtt{25134}\ \mathtt{31254}\ \mathtt{34215}\ \mathtt{41325}$\\
27&$\mathtt{12345}\ \mathtt{13524}\ \mathtt{15342}\ \mathtt{21453}\ \mathtt{23154}\ \mathtt{25413}\ \mathtt{32415}\ \mathtt{32514}\ \mathtt{42135}\ \mathtt{43125}$\\
28&$\mathtt{12345}\ \mathtt{13542}\ \mathtt{14253}\ \mathtt{15243}\ \mathtt{21435}\ \mathtt{21534}\ \mathtt{23514}\ \mathtt{25413}\ \mathtt{31254}\ \mathtt{42315}$\\
29&$\mathtt{12345}\ \mathtt{13542}\ \mathtt{14253}\ \mathtt{15243}\ \mathtt{21534}\ \mathtt{25143}\ \mathtt{31254}\ \mathtt{32145}\ \mathtt{41325}\ \mathtt{42315}$\\
30&$\mathtt{12345}\ \mathtt{13542}\ \mathtt{14253}\ \mathtt{15324}\ \mathtt{21435}\ \mathtt{23154}\ \mathtt{25413}\ \mathtt{32514}\ \mathtt{34215}\ \mathtt{43125}$\\
31&$\mathtt{12345}\ \mathtt{13542}\ \mathtt{14325}\ \mathtt{14523}\ \mathtt{15324}\ \mathtt{23154}\ \mathtt{24135}\ \mathtt{25134}\ \mathtt{34215}\ \mathtt{35214}$\\
32&$\mathtt{12345}\ \mathtt{13542}\ \mathtt{14325}\ \mathtt{15243}\ \mathtt{21534}\ \mathtt{23514}\ \mathtt{24135}\ \mathtt{31254}\ \mathtt{32145}\ \mathtt{42315}$\\
33&$\mathtt{12345}\ \mathtt{13542}\ \mathtt{14325}\ \mathtt{15243}\ \mathtt{21534}\ \mathtt{25413}\ \mathtt{31254}\ \mathtt{32415}\ \mathtt{35214}\ \mathtt{42315}$\\
34&$\mathtt{12345}\ \mathtt{13542}\ \mathtt{14523}\ \mathtt{15324}\ \mathtt{23154}\ \mathtt{24135}\ \mathtt{25143}\ \mathtt{32145}\ \mathtt{35124}\ \mathtt{43125}$\\
35&$\mathtt{12345}\ \mathtt{13542}\ \mathtt{14532}\ \mathtt{15243}\ \mathtt{21435}\ \mathtt{21534}\ \mathtt{31254}\ \mathtt{32415}\ \mathtt{41325}\ \mathtt{42315}$\\
36&$\mathtt{12345}\ \mathtt{13542}\ \mathtt{14532}\ \mathtt{15243}\ \mathtt{21534}\ \mathtt{24135}\ \mathtt{25143}\ \mathtt{31254}\ \mathtt{35214}\ \mathtt{42315}$\\
37&$\mathtt{12345}\ \mathtt{14253}\ \mathtt{14352}\ \mathtt{15324}\ \mathtt{21354}\ \mathtt{24315}\ \mathtt{24513}\ \mathtt{32514}\ \mathtt{34215}\ \mathtt{41325}$\\
38&$\mathtt{12345}\ \mathtt{14253}\ \mathtt{14352}\ \mathtt{15324}\ \mathtt{21354}\ \mathtt{24315}\ \mathtt{25143}\ \mathtt{31254}\ \mathtt{32145}\ \mathtt{41325}$\\
39&$\mathtt{12345}\ \mathtt{14325}\ \mathtt{14523}\ \mathtt{15324}\ \mathtt{21354}\ \mathtt{24315}\ \mathtt{24513}\ \mathtt{25314}\ \mathtt{34215}\ \mathtt{35214}$\\
40&$\mathtt{12345}\ \mathtt{14325}\ \mathtt{15243}\ \mathtt{15342}\ \mathtt{21354}\ \mathtt{24153}\ \mathtt{25314}\ \mathtt{31254}\ \mathtt{32145}\ \mathtt{42315}$\\
41&$\mathtt{12345}\ \mathtt{14325}\ \mathtt{15243}\ \mathtt{15342}\ \mathtt{21453}\ \mathtt{23154}\ \mathtt{25413}\ \mathtt{32415}\ \mathtt{35214}\ \mathtt{42135}$\\
42&$\mathtt{12345}\ \mathtt{14325}\ \mathtt{15324}\ \mathtt{15423}\ \mathtt{21354}\ \mathtt{24315}\ \mathtt{25314}\ \mathtt{25413}\ \mathtt{34125}\ \mathtt{35124}$\\
43&$\mathtt{12345}\ \mathtt{14325}\ \mathtt{15324}\ \mathtt{15423}\ \mathtt{21453}\ \mathtt{24135}\ \mathtt{25134}\ \mathtt{31254}\ \mathtt{34215}\ \mathtt{35214}$\\
44&$\mathtt{12345}\ \mathtt{14352}\ \mathtt{15243}\ \mathtt{15342}\ \mathtt{21354}\ \mathtt{21453}\ \mathtt{31254}\ \mathtt{32415}\ \mathtt{41325}\ \mathtt{42315}$\\
45&$\mathtt{12345}\ \mathtt{14352}\ \mathtt{15324}\ \mathtt{15423}\ \mathtt{21354}\ \mathtt{21453}\ \mathtt{24315}\ \mathtt{31254}\ \mathtt{34215}\ \mathtt{41325}$\\
46&$\mathtt{12345}\ \mathtt{14352}\ \mathtt{15324}\ \mathtt{15423}\ \mathtt{21354}\ \mathtt{24315}\ \mathtt{25143}\ \mathtt{31425}\ \mathtt{35124}\ \mathtt{41325}$\\
\bottomrule
\end{longtable}
\endgroup

\section{Adding or deleting several foldover pairs}\label{sec:mod}

Here a foldover OofA-OA is given, and the designs obtained from it by adding or deleting foldover pairs are compared. Every addition of one pair gives the same determinant, and so does every deletion of one pair, by~(4.2) and~(4.3) of the paper. For several pairs the result depends on the inner products $\kappa$ of~(2.5) between them, as follows.

\begin{proposition}\label{thm:augment}
Let $D_0=H_0\cup\rev H_0$ be a foldover OofA-OA with $|H_0|=h_0$. Add $r$ foldover pairs, which may repeat pairs of $D_0$ or each other, let $\mathcal R=(a^{(1)},\ldots,a^{(r)})$ contain one order of each added pair, and put $K_{\mathcal R}=Z_{\mathcal R}\Sigma^{-1}Z_{\mathcal R}^\top$. Its $(i,j)$ entry is $\kappa(a^{(i)},a^{(j)})$ and its diagonal entries are $q$. The final number of pairs is $h=h_0+r$. For the resulting design,
\begin{equation}\label{eq:augment}
 \frac{\det M_{D_+}}{\det M_F}
 =\left(\frac{h_0}{h_0+r}\right)^q\det(I_r+K_{\mathcal R}/h_0)
 \leq\left(\frac{h_0}{h_0+r}\right)^q(1+q/h_0)^r,
\end{equation}
with equality if and only if $\kappa(a^{(i)},a^{(j)})=0$ for all $i\neq j$, and such an addition is $D$-optimal among all additions of $r$ pairs to $D_0$. If instead $r$ pairs are deleted, the final number is $h=h_0-r$. If $h_0\geq q$ and $h_0-r>0$, then $I_r-K_{\mathcal R}/h_0$ is nonnegative definite and
\begin{equation}\label{eq:delete}
 \frac{\det M_{D_-}}{\det M_F}
 =\left(\frac{h_0}{h_0-r}\right)^q\det(I_r-K_{\mathcal R}/h_0)
 \leq\left(\frac{h_0}{h_0-r}\right)^q(1-q/h_0)^r,
\end{equation}
where $K_{\mathcal R}$ is formed from the deleted pairs, again with equality if and only if $\kappa(a^{(i)},a^{(j)})=0$ for all $i\neq j$.
\end{proposition}
These statements are the matrix determinant lemma and Hadamard's inequality applied to $h_0\Sigma\pm Z_{\mathcal R}^\top Z_{\mathcal R}$. For deletion, $K_{\mathcal R}/h_0$ is a principal submatrix of an orthogonal projection matrix, which gives the nonnegative definiteness.

For two pairs only $|\kappa(a,b)|$ matters, a smaller value giving a larger determinant, and for more pairs the $r\times r$ matrix $K_{\mathcal R}$ takes the place of the $q\times q$ determinant. Orders with $\kappa(a,b)=0$ need not exist. If $m\equiv2\pmod4$ then $q$ and every $z(a)^\top z(b)$ are odd while $s(a)^\top s(b)$ is even, so $(m+1)z(a)^\top z(b)-s(a)^\top s(b)$ is an odd integer and~(2.5) gives $|\kappa(a,b)|\geq3/(m+1)$.

\begin{corollary}\label{cor:two_pair_arithmetic}
Under the hypotheses of Proposition~\ref{thm:augment}, if $m\equiv2\pmod4$ and $r=2$, then
\begin{equation}\label{eq:arith}
 \frac{\det M_{D_\pm}}{\det M_F}
 \leq\left(\frac{h_0}{h_0\pm2}\right)^q
 \left\{\left(1\pm\frac{q}{h_0}\right)^2-\frac9{(m+1)^2h_0^2}\right\},
\end{equation}
with equality precisely when the orders $a$ and $b$ of the two pairs satisfy $|\kappa(a,b)|=3/(m+1)$. This follows by substituting the preceding bound into the two-by-two determinant.
\end{corollary}

For five components and the array generated by~(3.1) of the paper, the foldover pairs of $12345$, $13542$ and $31254$ are not in the array and have $K_{\mathcal R}=10I_3$. Adding them is optimal among additions of three pairs and gives thirty runs with determinant ratio $(12/15)^{10}(22/12)^3$ and $\Deff \approx 0.963148$. No four pairs outside the array have $\kappa(a,b)=0$ for all their orders. Theorem~6 does not apply when $h=16$, and evaluating Proposition~\ref{thm:augment} over all $\binom{48}4=194\,580$ additions of four pairs outside the array shows that the largest value of $\det(I_4+K_{\mathcal R}/12)$ is $605/54\approx11.203704$, against $(22/12)^4\approx11.2971$. It is attained by the pairs of $\mathtt{12345}$, $\mathtt{13542}$, $\mathtt{21534}$ and $\mathtt{31254}$, whose matrix $K_{\mathcal R}$ has two symmetric off-diagonal entries of absolute value two and zeros elsewhere off the diagonal. Exactly 174 additions of four pairs attain this value, and the resulting 32-run designs have $\Deff\approx0.958994$. Since $K_{\mathcal R}$ has integer entries here, maximising $\det(12I_4+K_{\mathcal R})$ is an integer calculation, with maximum $232\,320$, which gives $232\,320/12^4=605/54$.

For the 48-run foldover OofA-OA with six components of Section~\ref{sec:designs}, the orders $123465$ and $125463$ are not in the array and have Kendall tau distance five and $S=18$, so $\kappa=15-30+108/7=3/7$. Adding their foldover pairs attains~\eqref{eq:arith}, with $\Deff \approx 0.985743$. Deleting the pairs of $123456$ and $153426$, also with $|\kappa|=3/7$, attains the bound for deletion with $\Deff \approx 0.959664$. There are $6!/2=360$ foldover pairs, of which 336 are not in the array, hence $\binom{336}2=56\,280$ additions of two pairs and $\binom{24}2=276$ deletions of two pairs. Among these, $12\,097$ additions and $85$ deletions have $|\kappa|=3/7$ and attain the bounds of Corollary~\ref{cor:two_pair_arithmetic}, and their determinants were verified in exact arithmetic. Designs that differ by a relabelling of the components are counted separately.

For the example with seven components in Section~5, the 24 designs obtained by deleting one pair from the 48-run array all have $\det(Z_H^\top Z_H)=2^{78}$, as stated in~(5.3). None of the $24\times23\times2497=1\,378\,344$ exchanges of a single foldover pair increases the determinant of such a design. Each replacement is one of the $2520-23$ foldover pairs that are not in the current design, including the deleted pair, so these counts concern designs without repeated pairs. There are 1248 exchanges that leave the determinant unchanged, so the calculation does not exclude sequences of exchanges that pass through such designs. The design with $\det(Z_H^\top Z_H)=155\cdot2^{71}$ has no foldover pair in common with the 48-run array, and $24\Sigma-Z_H^\top Z_H$ has rank three, so it is not a foldover OofA-OA with one pair deleted. It has $V=501$, compared with the minimum $V=420$ allowed by Lemma~3 at $h=23$ and attained by every deletion of one pair.

These comparisons use an integer formula and not rounded determinants. Let $u$ be the order of the deleted pair, $v$ the order of the pair to be exchanged and $w$ the order of a candidate replacement, and write $z_u=z(u)$, $z_v=z(v)$ and $z_w=z(w)$ for their PWO vectors. Put $a=8\kappa(u,v)$, $b=8\kappa(u,w)$ and $c=8\kappa(v,w)$, all integers. The determinant lemma applied to $24\Sigma-z_uz_u^\top-z_vz_v^\top+z_wz_w^\top$ gives the ratio of the determinant to that of the design with one pair deleted as
\[
 \frac{405\,8^3-360a^2+24(b^2+c^2)+2abc}{24^2\,3\,8^3}.
\]
Every numerator is compared exactly with the denominator.

\section{Designs attaining the bounds by ranking orthogonal arrays}\label{sec:ranking}

Further designs attaining the bounds of Lemma~1 are obtained by ranking the rows of an ordinary orthogonal array OA$(N,m,v,4)$. This $N\times m$ array has $v$ ordered symbols, and in every four columns each combination of symbols occurs $N/v^4$ times. Its strength refers to the columns of symbols and not to PWO columns.

\begin{corollary}\label{cor:ranking}
Let $m\geq4$ and $v\geq2$, and suppose an OA$(N,m,v,4)$ exists. Fix an order $a$ of the components. In each row of the array, arrange the components in increasing order of their symbols, components with the same symbol being placed in the same relative order as in $a$. Add the reverse of every resulting order. The resulting design $D$ has $2N$ runs and
\begin{equation}\label{eq:rankingequality}
 G_D=(1-v^{-2})\Sigma+v^{-2}z(a)z(a)^\top,\qquad
 \theta(D)=v^{-2},
\end{equation}
so it attains the bounds (i)--(iii) of Lemma~1 at this value of $\theta$, with $\Deff(D)=\{(1-v^{-2})^{q-1}(1+(q-1)v^{-2})\}^{1/p}$. Orders that the construction produces more than once are kept as repeated runs.
\end{corollary}
\begin{proof}
Put $w=z(a)$. In a projection onto at most four components, the array symbols $U_i$ behave as independent uniform variables on $\{0,\ldots,v-1\}$, and a ranked comparison is $z_{ij}=\operatorname{sgn}(U_j-U_i)+\mathbf1_{\{U_i=U_j\}}w_{ij}$. For disjoint comparisons, independence gives second moment $v^{-2}w_{ij}w_{k\ell}$. For comparisons $(i,j)$ and $(i,k)$, the strict-comparison product has mean $v^{-1}\sum_{u=0}^{v-1}\{(v-1-2u)/v\}^2=(1-v^{-2})/3$, the mixed terms have mean zero, and a tie of all three symbols has probability $v^{-2}$, so the second moment is $(1-v^{-2})/3+v^{-2}w_{ij}w_{ik}$. For reversed indices, define $z_{ji}=-z_{ij}$ and $w_{ji}=-w_{ij}$; these identities give the remaining cases in which two comparisons share a component. The diagonal moments equal one, and the reverses remove the first moments. These identities give~\eqref{eq:rankingequality}, and $w^\top\Sigma^{-1}w=q$ gives $\theta(D)$ and the determinant.
\end{proof}

For example, the binary rows $(u_1,u_2,u_3,u_4,u_1+u_2+u_3+u_4)$, with addition modulo two, form an OA$(16,5,2,4)$. Breaking ties by $12345$ and adding the reverses gives 32 runs with 28 distinct orders, $h=16$ and $\theta=1/4$. The eigenvalues of $C_D$ are $13/4$ once and $3/4$ nine times, which gives $\Deff \approx 0.879655$. The rows $00000$, $00011$ and $01111$, for instance, all produce the order $12345$. On each set of four components the foldover pair of $12345$ occurs five times in $H$ and every other pair once, as Proposition~1 requires. The argument with moments fails when ties are broken by an order that depends on the row. To illustrate this, 4000 assignments of an order for breaking ties to each of the 16 rows were generated by the linear congruential generator $x\leftarrow(1103515245\,x+12345)\bmod 2^{31}$ with seed $5$, the order for each row being the permutation of rank $(x\gg7)\bmod 120$ in lexicographic order. Of these assignments, 649 gave 32 distinct runs and none attained the bounds, their efficiencies ranging from $0$ to $0.8908$, with $\theta\approx0.158$ for the largest. Other assignments are not excluded. The design with a common order for ties attains the bound at its value of $\theta$ but is not optimal among 32-run designs, since the best addition of four pairs to the 24-run array in Section~\ref{sec:mod} has efficiency $0.9590$.

\citet{tsai2022} also builds designs from combinatorial arrays. He juxtaposes Latin squares, each consisting of the $m$ orders obtained from one order by adding a constant to every label modulo $m$, and searches among such designs for OofA-OAs. The designs obtained here are not OofA-OAs, since $G_D$ is $(1-v^{-2})\Sigma$ plus a matrix of rank one, with $v^{-2}$ fixed by the number of symbols.

\section{Enumerations for three and four components}\label{sec:smallcases}

The following complete enumerations support the remarks in Sections~4.4 and~5 of the paper and give an independent check of Theorem~1 in its smallest cases.

For four components and $n=13$ runs, Corollary~1 gives $\det(X_D^\top X_D)\leq(n-1)^{p-1}(n+p-1)\det M_F=12^6\cdot19\cdot125/729=9\,728\,000$. Over all $\binom{24}{13}=2\,496\,144$ designs with distinct runs the maximum is exactly $9\,728\,000$, attained by 240 designs. Exactly 20 of the sets of twelve orders are OofA-OAs, and each of the 240 designs is one of these arrays with one of the twelve orders outside it added. With repeated runs allowed, a design attaining the bound must by Proposition~1, applied to $D\cup\rev D$, contain one run from each foldover pair except two from the pair of $a$. Enumerating all $24\times3\times2^{11}$ designs of this form gives, for each $a$, ten designs that are an array not containing $a$ with $a$ added, ten that are an array containing $a$ with $a$ repeated, and one that contains $\rev a$ twice and $a$ not at all. The last design contains no OofA-OA among its subsets of twelve runs, because removing a run $r$ other than $\rev a$ leaves the matrix $12\Sigma+z(a)z(a)^\top-z(r)z(r)^\top\neq12\Sigma$, and removing one copy of $\rev a$ leaves PWO column sums $2z(a)$. There are thus 24 designs in all that attain the bound and contain neither the order $a$ of their moment matrix nor any OofA-OA. An example with $a=1234$ is
\[
\mathtt{4321}\ \mathtt{4321}\ \mathtt{1243}\ \mathtt{1324}\ \mathtt{1342}\ \mathtt{1423}\ \mathtt{2341}\ \mathtt{2134}\ \mathtt{3412}\ \mathtt{2314}\ \mathtt{2413}\ \mathtt{3124}\ \mathtt{4123},
\]
with $X_D^\top X_D=12M_F+x(1234)x(1234)^\top$. Exchange searches with random restarts and repeated runs allowed never exceeded the bound.

For four components and $h=13$ foldover pairs, Theorem~6 gives $\det(Z_H^\top Z_H)\leq(h-1)^{q-1}(h+q-1)\det\Sigma=12^5\cdot18\cdot125/729=768\,000$. Over all $\binom{24}{13}$ multisets of 13 pairs from the twelve foldover pairs the maximum is exactly $768\,000$, attained by the twelve multisets consisting of the full design with one pair repeated.

\begin{table}[tb]
\centering
\caption{Designs with four components and $n$ distinct runs. The table gives the maximum of $\det(X_D^\top X_D)$, its efficiency, the number of designs attaining it, and how many of them are a 12-run OofA-OA with $n-12$ runs added or $12-n$ runs deleted. The bound is that of Corollary~1, which applies for every odd $n$.}
\label{tab:m4ordinary}
\begin{tabular}{rrrrlr}
\toprule
$n$ & $\max\det(X_D^\top X_D)$ & $\Deff$ & optima & array $\pm$ runs & bound\\
\midrule
7 & 65536 & 0.8961 & 24 & 0 of 24 & 0.9572\\
8 & 172032 & 0.9000 & 120 & 0 of 120 & --\\
9 & 442368 & 0.9156 & 320 & 160 of 320 & 0.9724\\
10 & 1081344 & 0.9363 & 360 & 0 of 360 & --\\
11 & 2560000 & 0.9627 & 240 & 240 of 240 & 0.9807\\
12 & 6144000 & 1.0000 & 20 & 20 of 20 & --\\
13 & 9728000 & 0.9857 & 240 & 240 of 240 & \textbf{0.9857}\\
14 & 15400960 & 0.9774 & 360 & 360 of 360 & --\\
15 & 24379392 & 0.9741 & 160 & 160 of 160 & 0.9889\\
16 & 38477824 & 0.9747 & 120 & 0 of 120 & --\\
17 & 61456384 & 0.9809 & 24 & 0 of 24 & 0.9913\\
18 & 91193344 & 0.9801 & 24 & 0 of 24 & --\\
19 & 134217728 & 0.9812 & 120 & 120 of 120 & 0.9929\\
20 & 194510848 & 0.9829 & 180 & 180 of 180 & --\\
21 & 279379968 & 0.9858 & 40 & 40 of 40 & 0.9941\\
22 & 394526720 & 0.9885 & 60 & 60 of 60 & --\\
23 & 557056000 & 0.9933 & 24 & 24 of 24 & 0.9950\\
24 & 786432000 & 1.0000 & 1 & 1 of 1 & --\\
\bottomrule
\end{tabular}
\end{table}

Table~\ref{tab:m4ordinary} gives, for four components and every run size from 7 to 24, the maximum of $\det(X_D^\top X_D)$ over all designs with distinct runs, the number of designs attaining it, and whether they are obtained from one of the 20 OofA-OAs with 12 runs by adding or deleting runs. The maximum at 7 runs, where the design is saturated, was found by \citet{zhao2021} in a complete search, and the efficiencies at 7, 13 and 19 runs agree with those reported by \citet[Table~6.1]{winker2020} for designs found by threshold accepting. The optimal designs are arrays with runs added or deleted at every run size from 11 to 15 and from 19 to 24, but not at 16, 17 and 18 runs, nor below 11 runs apart from half of the optimal designs at 9 runs. At 14 runs the 360 optimal designs are an array with two added orders $a,b$ satisfying $\kappa(a,b)=-6/5$, the value that minimises $|1+\kappa(a,b)|$ among pairs of orders outside the array; the other values of $\kappa(a,b)$ give determinants $15\,360\,000$, $15\,319\,040$, $15\,196\,160$ and $14\,909\,440$. For 14 runs the restriction to distinct runs can be removed. A design with a repeated run can be relabelled so that this run is $1234$, which leaves $\binom{35}{12}=834\,451\,800$ choices of the other twelve runs, and over all of them the maximum of $\det(X_D^\top X_D)$ is again $15\,400\,960$. The best additions of two runs to an array are therefore $D$-optimal among all 14-run designs. Some of the designs with a repeated run that attain the maximum contain no OofA-OA among their subsets of twelve runs. An example is
\[
\begin{array}{lllllll}
\mathtt{1234}&\mathtt{1234}&\mathtt{1432}&\mathtt{2143}&\mathtt{2431}&\mathtt{3124}&\mathtt{3142}\\
\mathtt{3214}&\mathtt{3241}&\mathtt{3421}&\mathtt{4132}&\mathtt{4213}&\mathtt{4231}&\mathtt{4312}
\end{array}
\]
for which $X_D^\top X_D=12M_F+x(3124)x(3124)^\top+x(4321)x(4321)^\top$. At 11 runs the 240 optimal designs are arrays with one run deleted, and at 23 runs the 24 optimal designs are the full design with one run deleted. For 11 runs the restriction to distinct runs can be removed. After one run is fixed as $1234$ by relabelling the components, there are $\binom{33}{10}=92\,561\,040$ designs with repeated runs allowed, and over all of them the maximum of $\det(X_D^\top X_D)$ is again $2\,560\,000$, attained only by the arrays with one run deleted. For 11 runs, then, an array with one run deleted is $D$-optimal among all designs, as Theorem~2 of \citet{tsai2022} states, and for 23 runs the full design with one run deleted is $D$-optimal among designs with distinct runs. The bound of Corollary~1 is attained at 13 runs and not attained at 7, 9, 11, 15, 17, 19, 21 and 23 runs. Over all multisets of 14 foldover pairs the maximum of $\det(Z_H^\top Z_H)$ is $1\,152\,000$, attained by the six multisets consisting of the full design and two added pairs whose orders have $\kappa(a,b)=0$, in agreement with Proposition~\ref{thm:augment}.

Whether the bounds are attained at a given run size and what the smallest such run size is are separate questions. For nine components the 49-run design of Table~2 is the smallest design attaining the bounds of Corollary~1. A design attaining them has $X_D^\top X_D=(n-1)M_F+x(a)x(a)^\top$, so Proposition~1 applied to $D\cup\rev D$ gives $n\equiv1\pmod{12}$, and the only candidate below 49 with $n\geq p=37$ is $n=37$. Then $X_D$ is square and $\det(X_D^\top X_D)=(\det X_D)^2$ would be a perfect square, whereas the required matrix has $\det(X_D^\top X_D)=36^{36}\cdot73\cdot10^8/3^{36}=73\,(12^{18}10^4)^2$, which is not. For ten components the congruence allows a design with 49 runs attaining the bounds, and we do not know whether one exists. For foldover designs the conditions $h\geq q$ and $h\equiv1\pmod{12}$ give at least 74 runs for nine components and 98 for ten, so the smallest foldover design with nine components attaining the bounds of Theorem~6 has 74 or 98 runs. The argument with the perfect square does not exclude $h=37$, because $Z_H$ is then $37\times36$ and not square. Estimability alone requires only $p$ runs without the foldover restriction, or $2q$ runs with it, as in Proposition~2.

For three components, Lemma~3 uses only the contribution of the four-element set, since every two of the four triples on $\{0,1,2,3\}$ share two labels, and the bounds of Theorem~6 and Corollary~1 are attained at $n=7$, $n=13$ and at $h=4,7,10$. Over all multisets of $n$ orders the maximum of $\det(X_D^\top X_D)$ is $1280$ at $n=7$ and $16\,384$ at $n=13$, each attained by the six designs consisting of copies of the full design and one further run. Over all multisets of $h$ pairs the maximum of $\det(Z_H^\top Z_H)$ is $32$, $192$ and $576$ at $h=4,7,10$, attained in each case by the three designs consisting of copies of the full design and one further pair. At $h=2$ and $h=8$ the bound is not attained, as the congruence $h\equiv1\pmod3$ requires. When one run is deleted from one, two or three copies of the full design, $\det(X_D^\top X_D)$ is $256$, $8192$ and $48\,384$, with $n=5$, $11$ and $17$. Over all multisets of $n$ orders these are the maxima, and each is attained only by the six designs of this form.

\section{Proof of Theorem~7 and its limits}\label{sec:pairingproofs}
We use the notation of Section~6.1 of the paper. There are finitely many design points, arranged in pairs $\{a_i,\rho a_i\}$ with regression vectors $f_i$ and $-f_i$, the numbers of runs at the two points of pair $i$ are $r_i$ and $s_i$ with total $t_i=r_i+s_i$, and $J$ is the matrix in~(6.3). All designs that are compared have the same total number of runs $n$.

\begin{lemma}\label{lem:oddpairrounding}
Suppose that $J$ is positive definite and that two pairs have odd totals $t_i,t_j$, with $|r_i-s_i|=|r_j-s_j|=1$. Consider the two designs with the same number of runs obtained by replacing these totals by $(t_i+1,t_j-1)$ and $(t_i-1,t_j+1)$ and dividing the runs equally within both pairs. If their information matrices are $J_+$ and $J_-$, then both are positive definite and
\begin{equation}\label{eq:roundingcovariance}
 \frac{J_+^{-1}+J_-^{-1}}2\preceq J^{-1}.
\end{equation}
Moreover,
\begin{equation}\label{eq:strictlogdetrounding}
 \frac{\log\det J_++\log\det J_-}{2}>\log\det J.
\end{equation}
\end{lemma}

\begin{proof}
Let $f_i,f_j$ be the two regression vectors and let $J_0$ contain the
contributions from the other pairs. Put
\[
 A=J_0+t_if_if_i^\top+t_jf_jf_j^\top,\qquad
 \Theta=f_if_i^\top-f_jf_j^\top,\qquad
 \Gamma=t_i^{-1}f_if_i^\top+t_j^{-1}f_jf_j^\top.
\]
Then $J=A-\Gamma$ and $J_\pm=A\pm\Theta$. The positive weights in $J$ remain
positive in both $J_+$ and $J_-$: an odd total at least three is replaced
by an even total at least two, and a total of one had weight zero in $J$.
Thus positive definiteness is preserved.

Write $\Phi_{ij}=[f_i\ f_j]$ and $T_{ij}=\operatorname{diag}(t_i,t_j)$. Since
$A\succeq\Phi_{ij}T_{ij}\Phi_{ij}^\top$, we have $\Phi_{ij}^\top A^{-1}\Phi_{ij}\preceq T_{ij}^{-1}$.
Indeed, $A^{-1/2}\Phi_{ij}T_{ij}^{1/2}$ has all singular values at most one. With
$S=\operatorname{diag}(1,-1)$, so that $\Theta=\Phi_{ij}S\Phi_{ij}^\top$, this implies
\[
 \Theta A^{-1}\Theta
 =\Phi_{ij}S(\Phi_{ij}^\top A^{-1}\Phi_{ij})S\Phi_{ij}^\top
 \preceq \Phi_{ij}ST_{ij}^{-1}S\Phi_{ij}^\top=\Gamma.
\]
The matrix harmonic-mean identity gives
\[
 \left\{\frac{(A+\Theta)^{-1}+(A-\Theta)^{-1}}2\right\}^{-1}
 =A-\Theta A^{-1}\Theta\succeq A-\Gamma=J.
\]
Inverting proves \eqref{eq:roundingcovariance}. The identity follows, for
example, by conjugating with $A^{-1/2}$ and using
$\{(I+\Xi)^{-1}+(I-\Xi)^{-1}\}/2=(I-\Xi^2)^{-1}$ for $\Xi=A^{-1/2}\Theta A^{-1/2}$.

Concavity and monotonicity of $\log\det$ now yield
\[
 \begin{aligned}
 \frac{\log\det J_++\log\det J_-}{2}
 &=-\frac{\log\det J_+^{-1}+\log\det J_-^{-1}}2
 \\ &\geq-\log\det\left\{\frac{J_+^{-1}+J_-^{-1}}2\right\}
 \geq\log\det J.
 \end{aligned}
\]
If $\Theta\ne0$, the first inequality is strict by strict concavity. If $\Theta=0$,
then $J_+=J_-=A=J+\Gamma$ and $\Gamma\ne0$, so the last comparison is strict.
This proves \eqref{eq:strictlogdetrounding}.
\end{proof}
\subsection{Proof of Theorem~7}

\begin{proof}
A design with $2q$ runs and positive definite $J$ exists: choose $q$ linearly independent vectors $f_i$ and run both points of each of these pairs once. Adding runs preserves positive definiteness, so a $D$-optimal design with $n\geq2q$ runs has a positive determinant. Dividing the runs within a pair as equally as possible increases $J$ in the Loewner order, and strictly increases its determinant if the two numbers differed by more than one. Thus $|r_i-s_i|\leq1$ for every pair of a $D$-optimal design.

The number of pairs with an odd total has the same parity as $n$. If at least two totals are odd, Lemma~\ref{lem:oddpairrounding} gives a design with the same number of runs and a strictly larger determinant. For even $n$ there can therefore be no odd total, and for odd $n$ there is exactly one. For odd $n$ this total cannot be one. A pair with a single run carries no information, and moving this run to either point of a pair with total $2k>0$ raises the coefficient of that pair from $2k$ to $2k+1-(2k+1)^{-1}$ and strictly increases the determinant. This proves (i) and (ii).

For (iii), let $\Lambda$ be any concave function of the covariance matrix that is nondecreasing in the Loewner order; $\tr(WJ^{-1})$ with a fixed nonnegative definite $W$ is such a function of $J^{-1}$. Dividing the runs equally within a pair cannot increase $\Lambda$. When two totals are odd, concavity and monotonicity give
\[
 \frac{\Lambda(J_+^{-1})+\Lambda(J_-^{-1})}{2}
 \leq\Lambda\left\{\frac{J_+^{-1}+J_-^{-1}}2\right\}
 \leq\Lambda(J^{-1}).
\]
so at least one of the two designs of Lemma~\ref{lem:oddpairrounding} has no larger value of $\Lambda$. Repeating the step removes two odd totals at a time. A remaining pair with a single run is treated as above without increasing $\Lambda$. Since there are finitely many designs with $n$ runs, an optimal design exists, and applying these steps to it gives an optimal design of the stated form.
\end{proof}
\subsection{Limits of Theorem~7}
Part (iii) cannot be strengthened to say that every design minimising $\tr(WJ^{-1})$ has the stated form, and the theorem does not extend to the criterion of the smallest eigenvalue. For three components the vectors $z(a)$ of the three foldover pairs are linearly independent and, after multiplication by $\Sigma^{-1/2}$, mutually orthogonal with squared length three. Consequently the eigenvalues of $\Sigma^{-1/2}J\Sigma^{-1/2}$ are $3\omega(t_1),3\omega(t_2),3\omega(t_3)$ when the runs are divided as equally as possible within the pairs. At $n=8$ the only totals with positive definite $J$ are $(2,3,3)$ and $(2,2,4)$, up to permutation, and both give
\[
 \operatorname{tr}(\Sigma J^{-1})=\frac5{12}.
\]
Both are therefore optimal for this criterion, but only the second has $r_i=s_i$ for all $i$. At $n=10$, the totals $(3,3,4)$ give the eigenvalues $(8,8,12)$, whereas every design with $r_i=s_i$ for all $i$ and positive definite $J$ has smallest eigenvalue at most $6$. Thus no design of the form in Theorem~7(i) maximises the smallest eigenvalue. This example also shows that the inequality~\eqref{eq:roundingcovariance} for the average of the two covariance matrices does not imply that one of the two designs is better than the original design in the Loewner order.

\subsection{Reduction for an odd number of runs}
For $n=2h+1$, consider each positive definite $A(k)=\sum_i k_i f_if_i^\top$ with nonnegative integers $k_i$ summing to $h$. Theorem~7 and the matrix determinant lemma give
\begin{equation}\label{eq:oddpairedreduction}
 \max_{\sum_i(r_i+s_i)=2h+1}\det J
 =2^q\max_k\left[\det A(k)
 \max_{i:k_i>0}\left\{1+\frac{k_i}{2k_i+1}
 f_i^\top A(k)^{-1}f_i\right\}\right].
\end{equation}
For a given $A(k)$, the best run to replicate maximises the expression in braces, and the formula does not require all allocations of the $2h+1$ runs to be enumerated. For $\tr(WJ^{-1})$, with $A=A(k)$ and one further run on pair $i$, the Sherman--Morrison formula gives
\[
 \operatorname{tr}(WJ^{-1})
 =\frac12\left\{\operatorname{tr}(WA^{-1})
 -\frac{c_i f_i^\top A^{-1}WA^{-1}f_i}
 {1+c_i f_i^\top A^{-1}f_i}\right\},\qquad
 c_i=\frac{k_i}{2k_i+1},
\]
so that the optimal design for odd $n$ is found by minimising this expression over $k$ and over $i$ with $k_i>0$.

\subsection{An example with an odd number of runs}
Take four pairs with $f_1=(3,0)^\top$, $f_2=(1,2)^\top$, $f_3=(4,1)^\top$ and $f_4=(4,2)^\top$. With $h=5$ there are $\binom83=56$ allocations $k$. Direct evaluation gives $\max\det A(k)=303$, attained at $k=(1,2,1,1)$. The maximum of the expression in square brackets in~\eqref{eq:oddpairedreduction} is $1768/5$, attained at $k=(0,2,2,1)$ by replicating a run of the second pair, and for this $k$ we have $\det A(k)=300$. The best replicate added to the design with $k=(1,2,1,1)$ gives only $1757/5$. Thus the $(n-1)$-run design contained in a $D$-optimal design with an odd number $n$ of runs need not be $D$-optimal for $n-1$ runs. This example is for the general model and does not use PWO factors.

\section{Proofs of Theorems~3, 4 and~5}\label{sec:residueproofs}
\subsection{Proof of Theorem~3}

\begin{proof}
Replacing an order by its reverse does not change $V(H)$, so we represent each of the twelve foldover pairs by its lexicographically smaller order. Arrange these twelve orders in six couples as follows, with an even permutation on the left and an odd one on the right of each couple:
\[
\begin{array}{c|c@{\qquad}c|c@{\qquad}c|c}
1234&2413&1342&3214&1423&2134\\
2143&1324&2314&1243&3124&1432.
\end{array}
\]
In the circular order $(0,a_1,a_2,a_3,a_4)$, the two orders of a couple correspond to the two complementary Hamilton cycles on the five labels. Let $c_j^+,c_j^-$ be the multiplicities of the two orders of couple $j$ and put $t_j=c_j^++c_j^-$ and $e_j=c_j^+-c_j^-$, for $j=1,\ldots,6$. Then $\sum_jt_j=h$, $e_j\equiv t_j\pmod2$ and $|e_j|\leq t_j$.

Let $c$ be the vector of the twelve multiplicities, $s$ the vector with $s_a=+1$ for even and $s_a=-1$ for odd permutations, and $P$ the permutation matrix that interchanges the two orders of each couple. Formula~(2.5) gives
\[
 25\bigl[\kappa(a,b)^2-6\bigr]_{a,b}
 =864I_{12}-60\mathbf1\mathbf1^\top-54ss^\top-144P.
\]
Indeed, the diagonal entries on the left are $750$, different orders of the same parity give $-114$, and orders of opposite parity give $-6$, except for the two orders of a couple, which give $-150$. These values follow by substituting $K$ and $S$ in~(2.5). Consequently
\begin{equation}\label{eq:m4quadratic}
25V(H)=360\sum_{j=1}^6t_j^2-60h^2
       +504\sum_{j=1}^6e_j^2-54\left(\sum_{j=1}^6e_j\right)^2.
\end{equation}

Let $\xi$ be the number of odd $t_j$.  The part involving $e$ can be written
\[
 180\sum_je_j^2+54\sum_{i<j}(e_i-e_j)^2
 \geq180\xi+54\xi(6-\xi)=504\xi-54\xi^2.
\]
The first inequality uses $e_j\equiv t_j\pmod2$: each odd coordinate
has square at least one, and each odd--even pair has squared difference
at least one.  Equality is attained by taking $e_j=0$ on the even
coordinates and the same sign, $1$ or $-1$, on all odd coordinates.

Write $h=6u+w$, where $0\leq w<6$, and let
$\xi_0=w$ when $u$ is even and $\xi_0=6-w$ when $u$ is odd.  For any integer
$k$, with $\epsilon(k)=\mathbf1_{\{k\text{ odd}\}}$,
$k(k-1)\geq|k-\epsilon(k)|$.  Applying this inequality to
$k_j=t_j-u$ and using $\sum_jk_j=w$ gives
\[
 \sum_jt_j^2\geq6u^2+(2u+1)w+|\xi-\xi_0|.
\]
This lower bound holds even if the nonnegativity restriction on $t$ is
removed.  Substitution in \eqref{eq:m4quadratic} yields
\begin{equation}\label{eq:m4finiteparity}
 25V(H)\geq60w(6-w)
 +\min_{\substack{0\leq \xi\leq6\\\xi\equiv h\ (\mathrm{mod}\ 2)}}
       \{360|\xi-\xi_0|+504\xi-54\xi^2\}.
\end{equation}
Evaluating the right-hand side for the possible values of $\xi$ gives precisely $25\mu_r$ in~(4.7). This is the evaluation of a proved inequality and not an enumeration of designs.

The bounds are attained as follows. For $r=1$, take any foldover pair. For $r=2$, take the two foldover pairs of one couple. For $r=3,4,5,6$, take any $r$ foldover pairs whose orders have the same parity. For $r=7,\ldots,11$, take the foldover pairs left out by one of these choices for $12-r$. For $r=0$ take none. These choices attain \eqref{eq:m4finiteparity}. The minimising value of $\xi$ in \eqref{eq:m4finiteparity} is unique for each $r$: for $r=0,\ldots,11$ it is $0,1,0,3,4,5,6,5,4,3,0,1$. Equality in the two inequalities above therefore forces exactly the stated multiplicities. More explicitly, except at $r=2,10$ the totals of the couples take only two consecutive integer values and are as equal as possible; at $r=2$ one total is $2k+2$ and the other five are $2k$, and at $r=10$ one is $2k$ and the other five are $2k+2$. The condition for equality on $e$ puts the extra run of every couple with an odd total on permutations of the same parity. This also gives the stated numbers of designs. Adding one copy of each of the twelve foldover pairs changes $h$ to $h+12$ and $Z_H^\top Z_H$ to $Z_H^\top Z_H+12\Sigma$, and so leaves $V$ unchanged. This proves the theorem for every $h$, with nonnegative multiplicities attaining the lower bound.
\end{proof}
\subsection{Proof of Theorem~4}

\begin{proof}
Use the orientations of triples of Lemma~2 on
the $\nu$ labels $\{0,1,\ldots,m\}$.  For $j=0,1,2$, let $S_j$ be
the sum of squared entries of $U^\top U-hR$ over ordered pairs of
distinct triples having $j$ common labels.  Then
\[
 V(H)=\frac9{\nu^2}(S_0+S_1+S_2).
\]
When $j\leq1$, the reference cross-product is zero.  The corresponding
empirical cross-product is a sum of $h$ signs and therefore has square
at least $\epsilon_h$.  There are $20\binom\nu6$ ordered pairs with
no common label and $30\binom\nu5$ with one common label, giving
\begin{equation}\label{eq:overlapparity}
 S_0\geq20\binom\nu6\epsilon_h,\qquad
 S_1\geq30\binom\nu5\epsilon_h.
\end{equation}

For each five-label set $W$, restrict every circular order to $W$,
choose any one of its labels as anchor, and read the other four labels
after that anchor.  This produces a multiset of four-component orders.
Changing the anchor or relabelling the five labels only changes fixed
sign conventions in the triple coordinates, and hence preserves the
squared sums.  Theorem~3 therefore
implies that the sum of the squared entries of $U^\top U-hR$ for pairs of triples
contained in $W$ is at least $25\mu_{h\bmod12}/9$.
Summing over $W$, a triple pair with one common label is counted once,
and a pair with two common labels is counted $\nu-4$ times.  Thus
\begin{equation}\label{eq:fivelabelsum}
 S_1+(\nu-4)S_2\geq\binom\nu5\frac{25}{9}\mu_{h\bmod12}.
\end{equation}

There is a further restriction if $h\equiv2\pmod4$ and $\nu\geq6$.
Take three triples that meet pairwise in one label, with the three
intersection labels distinct and the union containing six labels.
Write $y_{12},y_{13},y_{23}$ for their empirical cross-products.
For each observation the sum of their three sign products is either
$3$ or $-1$, so
\[
 y_{12}+y_{13}+y_{23}\equiv-h\equiv2\pmod4.
\]
Each $y_{ij}$ is even, and hence
$y_{12}^2+y_{13}^2+y_{23}^2\geq4$.
There are $120\binom\nu6$ such unordered sets of three triples:
choose the three intersection labels and then assign the other three
labels to the three triples.  Each unordered pair of triples sharing
one label belongs to $4(\nu-5)$ such sets.  Summing the inequality gives
\begin{equation}\label{eq:sixlabelparity}
 S_1\geq40\binom\nu5\qquad(h\equiv2\pmod4).
\end{equation}

Now
\[
 S_0+S_1+S_2
 =S_0+\frac{S_1+(\nu-4)S_2}{\nu-4}
       +\frac{\nu-5}{\nu-4}S_1.
\]
Use \eqref{eq:overlapparity} and \eqref{eq:fivelabelsum}, and substitute
\eqref{eq:sixlabelparity} for the lower bound on $S_1$ when it applies.
Simplification gives (4.8).  When $\nu=5$,
the last coefficient is zero and the bound for four components is
recovered.  The statements for particular values of $h$ modulo 12
follow by substituting (4.7); in particular the bound
is positive at $h\equiv6\pmod{12}$, where neither of the two
contributions in Lemma~3 is available.
\end{proof}
\subsection{Proof of Theorem~5}

\begin{proof}
Theorem~4 and (4.7) give $V(H)\geq q(q-1)$ for every odd $h$, with a strictly larger lower bound unless $h\equiv1,11\pmod{12}$. If equality holds, the orders induced on every five labels used in \eqref{eq:fivelabelsum} must attain the minimum for four components. In particular, this applies to each set $\{0\}\cup Q$, where $Q$ is a set of four components. By the description of the minimising designs in Theorem~3, the matrix $Z^\top Z$ of the orders induced on $Q$ is
\[
 A_Q=(h-\varepsilon)\Sigma_Q+
        \varepsilon w_Qw_Q^\top,
 \qquad
 \varepsilon=\begin{cases}1,&h\equiv1\pmod{12},\\-1,&h\equiv11\pmod{12},\end{cases}
\]
where $w_Q$ is the PWO vector of an order on $Q$.

Put $A=Z_H^\top Z_H$ and $F=\varepsilon\{A-(h-\varepsilon)\Sigma\}$.
The coordinates of $F$ are indexed by edges of the complete graph on
the $m$ components.  It has diagonal one, and every block indexed by
the six edges on four components is a rank-one sign matrix.  We show
that these local blocks determine a single global rank-one matrix.

For the star edges through component $1$, set $w_{12}=1$ and
$w_{1j}=F_{12,1j}$ for $j\geq3$.  Every three star edges have at most
four endpoints, so local rank one gives
$F_{1i,1j}=w_{1i}w_{1j}$.  For an edge $ij$ not incident to $1$, put
$w_{ij}=F_{ij,1i}w_{1i}$.  Applying local rank one to
$\{ij,1i,1k\}$, whose endpoints lie in at most four components, gives
$F_{ij,1k}=w_{ij}w_{1k}$ for every $k\ne1$.
For any two edges sharing a vertex, adjoin a star edge at that vertex;
the same argument gives $F_{ij,ik}=w_{ij}w_{ik}$.
For disjoint edges $ij$ and $k\ell$, use the connecting edge
$ik$.  All three lie on four components, so
\[
 F_{ij,k\ell}=F_{ij,ik}F_{ik,k\ell}=w_{ij}w_{k\ell}.
\]
Thus $F=ww^\top$ globally.

The averaged PWO triangle identity holds for $A$ with right-hand side
$h$, and for $\Sigma$ with right-hand side one.  Substituting
$A=(h-\varepsilon)\Sigma+\varepsilon ww^\top$ shows that
\[
 w_{ij}w_{ik}-w_{ij}w_{jk}+w_{ik}w_{jk}=1\qquad(i<j<k).
\]
Hence the signs of $w$ orient a tournament with no cyclic triangle,
which defines a total order $a$, and $w=z(a)$.
Conversely, either displayed information matrix has
$\Sigma^{-1/2}(A-h\Sigma)\Sigma^{-1/2}
=\varepsilon(vv^\top-I_q)$ with $v^\top v=q$, giving
$V(H)=q(q-1)$ directly.
\end{proof}
The following corollary of Theorem~4 is described in words in Section~4.3 of the paper.
\begin{corollary}\label{cor:allresiduecriteria}
For a foldover design with $m\geq4$ and $h>0$ foldover pairs, set
\[
 t_{m,h}=\frac1h\sqrt{\frac{B_m(h)}{q(q-1)}}.
\]
Then the bounds of Lemma~1(i)--(iii) hold with
$t_{m,h}$ in place of $\theta$ on their right-hand sides.  In particular,
\[
 \frac{\det M_D}{\det M_F}
 \leq(1-t_{m,h})^{q-1}\{1+(q-1)t_{m,h}\}.
\]
For an arbitrary design with $n$ runs, with $p=q+1$, put
\[
 \widetilde B_m(n)=B_m(n)+2q\epsilon_n,
 \qquad
 \widetilde t_{m,n}
 =\frac1n\sqrt{\frac{\widetilde B_m(n)}{p(p-1)}}.
\]
The same three bounds hold with $p$ in place of $q$ and
$\widetilde t_{m,n}$ on their right-hand sides.  Thus, for every run
size,
\[
 \frac{\det M_D}{\det M_F}
 \leq(1-\widetilde t_{m,n})^{p-1}
       \{1+(p-1)\widetilde t_{m,n}\}.
\]
If $t_{m,h}$ or $\widetilde t_{m,n}$ equals one, the lower bound on the trace of the inverse is taken as $+\infty$.
\end{corollary}
\subsection{Proof of Corollary~\ref{cor:allresiduecriteria}}

\begin{proof}
For foldover designs, $h^2\Delta=V(H)\geq B_m(h)$ and the three
right-hand sides are monotone in $\theta$.  For an arbitrary design the
identity used in the proof of Corollary~1 gives
\[
 n^2\tr\{(\mathcal C_D-I_p)^2\}
 =V(D)+2b_D^\top\Sigma^{-1}b_D.
\]
The representation by triples gives
$b_D^\top\Sigma^{-1}b_D\geq q\epsilon_n$, because every entry
of $Lb_D$ is a sum of $n$ signs.  Apply
Theorem~4 to $H=D$ and then
Lemma~1 with $p$ in place of $q$.
\end{proof}
The bound $B_m(h)$ need not be attained for $m>4$. When the condition on six labels is not used, equality requires that the orders induced on every five labels minimise $V$ for four components, and that the sums of products of orientations for triples sharing at most one label have absolute value one when $h$ is odd and zero when $h$ is even; the condition for triples sharing one label is needed only for $m>4$. When the condition on six labels is used, the requirement for triples sharing one label is replaced by equality in the sum of the inequalities leading to \eqref{eq:sixlabelparity}.

\section{Proof of Proposition~3 and two added runs for three components}\label{sec:tworunproofs}
\subsection{Proof of Proposition~3}
\begin{proof}
Put $u=M_F^{-1/2}x(a)$ and $v=M_F^{-1/2}x(b)$. Both squared norms equal $p$, and $u^\top v=c$. The two possibly nonzero eigenvalues of $uu^\top+vv^\top$ are $p+|c|$ and $p-|c|$. Hence the eigenvalues of $\mathcal C_D$ are $(n_0+p+|c|)/n$ and $(n_0+p-|c|)/n$, together with $n_0/n$ repeated $p-2$ times. These give (6.5)--(6.7), including the degenerate case $|c|=p$, and the statements on monotonicity follow directly. For $b=\rev a$ we have $c=1-q$. Two orders that differ by an interchange of adjacent components have $K=1$ and $S=2$, so~(2.5) gives $c(a,b)=q-5+12/(m+1)$, which equals one for $m=3$ and is positive for $m\geq4$. In every case it is strictly less than $q-1$, because $12/(m+1)<4$, which proves (ii). For (iii), (2.5) gives $(m+1)c(a,b)=(m+1)p+6\{S(a,b)-(m+1)K(a,b)\}$. The number $S(a,b)$ is even, because modulo two the sum of the squared differences of positions equals the sum of these differences, which is zero. When $m$ is odd, $m+1$ is even, so the expression in braces is even as well. This proves the congruence and the lower bound, and the values of $K$ and $S$ in Table~\ref{tab:tworunpairs} show that the bound is attained for $3\leq m\leq30$ without a search over pairs of orders. The statements on the criteria follow from (i).
\end{proof}
Table~\ref{tab:tworunpairs} lists, for each number of components from three to thirty, an order $b$ such that the identity order $a$ and $b$ attain the lower bound $|c(a,b)|=\delta_m$ of Proposition~3(iii), together with the values of $K(a,b)$, $S(a,b)$ and $c(a,b)$. Where the table gives block lengths, $b$ is obtained from $a$ by reversing blocks of consecutive components with these lengths, the first block listed being at the end of the order, so that the lengths $5,4,3,2$ for 14 components stand for $b=(2,1,5,4,3,9,8,7,6,14,13,12,11,10)$. Reversing a block of $t$ consecutive components adds $t(t-1)/2$ to $K$ and $t(t^2-1)/3$ to $S$. We do not know whether the bound is attained for every $m$.

\begin{table}[p]
\centering
\caption{Orders attaining the lower bound of Proposition~3(iii) for two added runs. The first order is $a=(1,\ldots,m)$ and the second is $b$; $K$ and $S$ are the Kendall tau distance and the sum of squared differences of positions between them. Block lengths are explained in the text.}
\label{tab:tworunpairs}
\begin{tabular}{rllrrl}
\toprule
$m$ & $b$ & $\delta_m$ & $K$ & $S$ & $c(a,b)$\\
\midrule
3 & $(1,3,2)$ & $1$ & 1 & 2 & $1$\\
4 & $(2,1,4,3)$ & $1/5$ & 2 & 4 & $-1/5$\\
5 & $(1,2,5,4,3)$ & $1$ & 3 & 8 & $1$\\
6 & $(1,4,6,3,5,2)$ & $2/7$ & 7 & 30 & $-2/7$\\
7 & $(1,2,5,7,6,4,3)$ & $1/2$ & 8 & 34 & $-1/2$\\
8 & $(1,2,3,7,8,6,5,4)$ & $1/3$ & 9 & 38 & $1/3$\\
9 & $(1,2,3,6,9,8,7,5,4)$ & $1/5$ & 12 & 58 & $-1/5$\\
10 & $(1,2,3,7,10,8,6,9,5,4)$ & $2/11$ & 16 & 92 & $2/11$\\
11 & $(1,2,3,4,7,11,10,9,8,6,5)$ & $0$ & 17 & 92 & $0$\\
12 & $(1,2,3,4,11,12,9,10,7,5,8,6)$ & $1/13$ & 23 & 154 & $1/13$\\
13 & $(1,2,3,4,5,13,12,11,8,7,10,9,6)$ & $1/7$ & 24 & 152 & $1/7$\\
14 & blocks 5, 4, 3, 2 & $0$ & 20 & 70 & $0$\\
15 & $(1,2,3,4,5,6,14,15,13,11,12,10,8,7,9)$ & $1/4$ & 32 & 230 & $1/4$\\
16 & $(1,2,3,4,5,6,16,15,11,10,14,9,13,12,8,7)$ & $1/17$ & 37 & 286 & $-1/17$\\
17 & $(1,2,3,4,5,6,7,17,16,15,13,10,14,12,9,11,8)$ & $1/3$ & 40 & 310 & $1/3$\\
18 & blocks 5, 5, 4, 4 & $2/19$ & 32 & 120 & $-2/19$\\
19 & blocks 7, 6, 2, 2 & $1/5$ & 38 & 186 & $-1/5$\\
20 & blocks 7, 5, 5 & $1/7$ & 41 & 192 & $-1/7$\\
21 & $(1,\ldots,9,21,20,19,17,18,15,14,13,11,16,12,10)$ & $1/11$ & 60 & 546 & $-1/11$\\
22 & blocks 8, 7, 2, 2 & $2/23$ & 51 & 284 & $2/23$\\
23 & blocks 9, 7 & $0$ & 57 & 352 & $0$\\
24 & blocks 8, 8, 3, 2 & $1/25$ & 60 & 346 & $1/25$\\
25 & blocks 12, 4, 2 & $1/13$ & 73 & 594 & $1/13$\\
26 & blocks 11, 5, 4, 2, 2 & $0$ & 73 & 504 & $0$\\
27 & blocks 9, 6, 6, 4, 2 & $1/7$ & 73 & 402 & $1/7$\\
28 & blocks 11, 6, 4, 3, 3 & $1/29$ & 82 & 546 & $-1/29$\\
29 & blocks 10, 8, 5, 3 & $1/5$ & 86 & 546 & $1/5$\\
30 & blocks 12, 6, 5, 3, 2 & $2/31$ & 95 & 692 & $-2/31$\\
\bottomrule
\end{tabular}
\end{table}
\subsection{Two added runs for three components}
\label{sec:threecomponenttwo}

For three components the comparison can be extended from additions of two runs to all designs. In this case the design given by Proposition~3(iii) is optimal among all designs of its run size.

\begin{proposition}\label{prop:threecomponenttwo}
For three components and $n=6k+2$, with integer $k\geq1$, the largest
value of $\det(X_D^\top X_D)$ over all designs with $n$ runs is
\[
 \max_D\det(X_D^\top X_D)=64k^2(2k+1)(6k+5).
\]
The designs attaining it consist of $k$ copies of each of the six orders
and one further order from each of two different foldover pairs. There
are twelve such designs.
\end{proposition}
\begin{proof}
The PWO vectors of the three foldover pairs are mutually orthogonal in
the inner product defined by $\Sigma^{-1}$, each with squared norm three.
After an orthogonal change of scaled coordinates, the six model vectors
are therefore $(1,\pm\sqrt3 e_i^\top)^\top$, $i=1,2,3$.
Let the two multiplicities on pair $i$ have sum $t_i$ and difference $b_i$.
For $t_i>0$, a Schur complement and $\det M_F=16/27$ give
\[
 \det(X_D^\top X_D)
 =16\left\{nt_1t_2t_3-b_1^2t_2t_3-b_2^2t_1t_3-b_3^2t_1t_2\right\}.
\]
A missing pair makes the determinant zero. For fixed positive $t_i$, the
determinant is largest when $b_i=0$ for even $t_i$ and $|b_i|=1$ for odd
$t_i$. Since $n$ is even, either all three totals are even or precisely two
are odd.

If all are even, their product is maximised at the totals
$(2k,2k,2k+2)$ up to permutation. The expression in braces is then
$n(2k)^2(2k+2)$.
If precisely two totals are odd, write the even one as $a$ and the other
two as $b,c$. The expression in braces becomes
\[
 a\{nbc-(b+c)\}\leq
 f(a):=a\{n(n-a)^2/4-(n-a)\}.
\]
Here $a$ is a positive even integer with $a\leq n-2$.
The derivative
\[
 f'(a)=\frac n4(n-a)(n-3a)-n+2a
\]
is positive for $a\leq2k$ and negative for $a\geq2k+2$ in this range.
Moreover $f(2k)-f(2k+2)=2k(n+2)>0$. Hence the largest value on the even
integers is $f(2k)$, attained only for
$a=2k$ and $b=c=2k+1$. This value exceeds the best all-even value by
$4k^2>0$. Substitution yields
\[
 16f(2k)=64k^2(2k+1)(6k+5).
\]
The pair with the even total has multiplicities $(k,k)$, and each of the
two pairs with odd totals has multiplicities $(k,k+1)$ in either order.
There are three choices for the first pair and two choices of the order
of the multiplicities for each of the other two, which gives the twelve
designs.
\end{proof}

\section{Further proofs}\label{sec:additionalproofs}
\subsection{Proof of Theorem~9}

\begin{proof}
Represent each of the sixty foldover pairs by its lexicographically smaller order. Relabelling the components permutes the PWO columns and changes the signs of some of them, and so preserves $|\det Z_H|$. Any selected order can therefore be relabelled to $12345$, and the other orders can then be replaced by their reverses independently, so it suffices to consider the $\binom{59}9=12\,565\,671\,261$ selections of ten orders containing $12345$. The branch-and-bound search in exact integer arithmetic described in Section~\ref{sec:enum} discards a branch only when its bound is strictly below $8192^2$, accounts for all selections, and finds exactly 46 selections with $|\det Z_H|=8192$ and none with a larger value, each verified by integer elimination. Hence $\max_D\det M_D=8192^2/10^{10}=2^{16}/5^{10}$ by~(2.2), and~(2.1) gives the efficiency.

For the classification, every permutation in $S_5$ is applied to each retained selection, each order is replaced by the smaller order of its foldover pair, and the sorted ten orders are compared. The smallest image is the canonical form, and this gives the five designs. If $a_i$ is the order of the stabiliser of design $i$ in $S_5$ and $c_i$ the number of retained selections equivalent to it, the orbit of design $i$ has $120/a_i$ designs. Since $S_5$ is transitive on the sixty foldover pairs and each design contains ten of them, one sixth of each orbit contains the pair of $12345$, so $120/a_i=6c_i$, and the sizes of the orbits add up to $276=6\times46$.

For each of the five designs put $W_H=(6I_q-B^\top B)Z_H^\top Z_H$, whose eigenvalues are twenty times those of $C_D$. Integer multiplication gives $(W_H-8I_q)(W_H-24I_q)(W_H-32I_q)=0$ for designs 1--3 and $(W_H-12I_q)(W_H-32I_q)=0$ for designs 4 and 5. The traces of the corresponding projection matrices give the multiplicities, the eigenvalues give the criteria, and counting gives the frequencies of the positions.
\end{proof}
\subsection{Proof of Proposition~1}

\begin{proof}
Extend the notation by $z_{ji}=-z_{ij}$. For distinct $i,j,k,\ell$, the indicator of the relative order $i\prec j\prec k\prec\ell$ is $(1+z_{ij})(1+z_{jk})(1+z_{k\ell})/8$. Adding the indicator of the reverse order cancels the odd terms, so the indicator of their foldover pair is
\[
 \frac{1+z_{ij}z_{jk}+z_{ij}z_{k\ell}+z_{jk}z_{k\ell}}4.
\]
Its expectation is an affine function of the entries of $G_D$. This expectation equals $1/12$ under $G_D=\Sigma$ and equals $\mathbf1_{\{b=[a|_Q]\}}$ under $G_D=z(a)z(a)^\top$. Substitution of~(4.11) and multiplication by $h$ prove~(4.12), for either sign of $\varepsilon$.

Conversely, the twelve counts on $Q$ determine the distribution of the relative order on $Q$ up to reversal, hence every product moment $z_{ij}z_{k\ell}$ with all four indices in $Q$, and every second PWO moment involves at most four components; the hypothesis $m\geq4$ is needed only so that every pair of PWO columns lies in some four-set. The counts~(4.12) therefore determine all entries of $G_D$ and give~(4.11). On a fixed $Q$, the distinguished count differs from every other count by $\varepsilon t$, so $t\in\mathbb Z$, and the other counts give $h-\varepsilon t\in12\mathbb Z$; for $\varepsilon=-1$ the distinguished count $(h-11t)/12$ must be nonnegative.
\end{proof}
\subsection{Proof of Theorem~2}

\begin{proof}
Equality in (i) with $0<\theta<1$ forces $C_D=(1-\theta)I_q+\theta vv^\top$ with $v^\top v=q$, by the case of equality in Lemma~1. Writing $u=\Sigma^{1/2}v$ gives $G_D=(1-\theta)\Sigma+\theta uu^\top$. Both $G_D$ and $\Sigma$ have unit diagonal, so every entry of $u$ is $1$ or $-1$. For $i<j<k$, the PWO identity of \citet{mee2020} is
\begin{equation}\label{eq:triangle}
 z_{ij}z_{ik}-z_{ij}z_{jk}+z_{ik}z_{jk}=1,
\end{equation}
whereas the two cyclic sign patterns give $-3$. Averaging~\eqref{eq:triangle} over the design and substituting the equality form, with $\Sigma_{ij,ik}=\Sigma_{ik,jk}=1/3$ and $\Sigma_{ij,jk}=-1/3$, shows that $u$ satisfies the same identity. The signs of $u$ orient one edge between each pair of components, and an orientation without cyclic triangles is transitive and defines a total order, so $u=z(a)$ for an order $a$. This proves necessity in~(4.1), and the eigenvalues of the right-hand side prove sufficiency. The design measure that puts weight $1-\theta$ on the uniform measure over all orders and weight $\theta$ on the pair $\{a,\rev a\}$ attains every $\theta\in[0,1]$. Part (ii) follows in the same way from the one smaller eigenvalue, since the diagonal and triangle arguments apply to a negative coefficient as well.
\end{proof}
\subsection{Proof of Corollary~1}
\begin{proof}
Put $A=Z_D^\top Z_D$ and $b_D=Z_D^\top\mathbf1\in\R^q$. The trace of $\mathcal C_D$ is $p$, because $\tr(\Sigma^{-1}A)/n$ equals $\sum_{a\in D}\kappa(a,a)/n=q$. The matrix $M_F^{-1}(M_D-M_F)$ has the blocks $0$ and $b_D^\top/n$ in its first row and $\Sigma^{-1}b_D/n$ and $\Sigma^{-1}(A/n-\Sigma)$ in its second, so the trace of its square gives
\[
 n^2\tr\{(\mathcal C_D-I_p)^2\}
 =\tr[\{\Sigma^{-1}(A-n\Sigma)\}^2]+2b_D^\top\Sigma^{-1}b_D.
\]
The first term is at least $q(q-1)$ by Theorem~5 for $m\geq4$ and by Lemma~3 for $m=3$, applied to $H=D$. Each entry of $Lb_D=\sum_{a\in D}u(a)$ is a sum of $n$ signs and hence an odd integer, so by~(4.4)
\[
 b_D^\top\Sigma^{-1}b_D=\frac3{m+1}\|Lb_D\|^2
 \geq\frac{3|\mathcal T|}{m+1}=q,
\]
and the total is at least $p(p-1)$, that is $\theta_p\geq1/n$, where $\theta_p^2=\tr\{(\mathcal C_D-I_p)^2\}/\{p(p-1)\}$. Lemma~1 applied to $\mathcal C_D$, with $p$ in place of $q$ and $\theta_p$ in place of $\theta$, gives the four bounds. In the case of equality, $\mathcal C_D=(1-\theta_p)I_p+\theta_pvv^\top$ gives $M_D=(1-\theta_p)M_F+\theta_pww^\top$ with $w=M_F^{1/2}v$. The unit diagonals force $w\in\{1,-1\}^p$, and we choose the sign of $w$ so that its first entry is one. The argument with the identity~\eqref{eq:triangle} in the proof of Theorem~2, applied to the PWO block, then gives $w=(1,z(a)^\top)^\top$ for an order $a$, and $ww^\top=x(a)x(a)^\top$. The congruences follow from the PWO block exactly as in Theorem~6, with Proposition~1 applied to $D\cup\rev D$. Adding one run to an OofA-OA gives $X_D^\top X_D=(n-1)M_F+x(a)x(a)^\top$, for which the eigenvalues of $\mathcal C_D$ are $(n+p-1)/n$ once and $(n-1)/n$ with multiplicity $p-1$, and this attains all four bounds. The run size of an OofA-OA is a multiple of twelve when $m\geq4$ and of six when $m=3$, as recalled in Section~2 of the paper, so $n$ is then odd and not divisible by three.
\end{proof}
\subsection{Designs with a run or a foldover pair deleted}\label{sec:deletionlist}
This section gives the designs of Section~5 of the paper and the calculations behind Table~3. All determinants were obtained by integer elimination.

\noindent\textit{Five components.} Up to relabelling of the components and reversal of all orders there are two OofA-OAs with five components and 12 runs \citep[Table~2]{schoenmee2023}. One array from each class is
\[
\begin{array}{llllll}
\mathtt{13245}&\mathtt{13542}&\mathtt{14253}&\mathtt{15243}&\mathtt{23415}&\mathtt{23514}\\
\mathtt{25413}&\mathtt{43215}&\mathtt{43512}&\mathtt{45213}&\mathtt{53214}&\mathtt{53412}
\end{array}
\]
and
\[
\begin{array}{llllll}
\mathtt{12453}&\mathtt{13452}&\mathtt{21534}&\mathtt{23451}&\mathtt{31524}&\mathtt{32514}\\
\mathtt{41325}&\mathtt{42315}&\mathtt{45321}&\mathtt{51423}&\mathtt{52413}&\mathtt{53412}
\end{array}
\]
Both have $X^\top X=12M_F$, and deleting any run from either of them gives $|\det X_D|=36\,864=2^{12}\cdot3^2$, as~(5.1) requires. The design~(5.2) has $|\det X_D|=45\,056=2^{12}\cdot11$. For this design the eigenvalues of $\mathcal C_D$ are $4/11$ twice, $12/11$ six times and three further values, approximately $0.6021$, $1.4021$ and $1.7232$. Corollary~1 gives $\det(X_D^\top X_D)\leq10^{10}\cdot21\cdot16/729$ for 11 runs, that is $|\det X_D|\leq67\,890$, so the largest value of $|\det X_D|$ lies between $45\,056$ and $67\,890$. An exchange search from $60\,000$ random starting designs found no design with $|\det X_D|>45\,056$. It found 2072 designs with this value, which fall into nine classes under relabelling and reversal. Other searches had reached this value before. \citet[Table~6.1]{winker2020} report the efficiency $0.90267$ for five components and 11 runs, obtained by threshold accepting, and the 11-run design that \citet{wangwang2023} give for the $D$-criterion, which is not equivalent to~(5.2), also has $|\det X_D|=45\,056$.

\noindent\textit{Six components.} By~(5.1), a 24-run OofA-OA with six components and one run deleted has $\det(X_D^\top X_D)=7^5\cdot2^{48}$. This was confirmed for the array of \citet[Example~3]{tsai2022}. The 23 orders
\[
\begin{array}{llllll}
\mathtt{126345}&\mathtt{132546}&\mathtt{145632}&\mathtt{162453}&\mathtt{243165}&\mathtt{243561}\\
\mathtt{256341}&\mathtt{261543}&\mathtt{315642}&\mathtt{346152}&\mathtt{352164}&\mathtt{364251}\\
\mathtt{421536}&\mathtt{452631}&\mathtt{463251}&\mathtt{514236}&\mathtt{532416}&\mathtt{536124}\\
\mathtt{541623}&\mathtt{623154}&\mathtt{641352}&\mathtt{652134}&\mathtt{654312}&
\end{array}
\]
have $\det(X_D^\top X_D)=3\cdot7^3\cdot263\cdot2^{44}$, which is larger by the factor $789/784$.

\noindent\textit{Seven components.} The foldover design with seven components and 46 runs of Section~5 of the paper has the half-design
\[
\begin{array}{llllll}
\mathtt{4267315}&\mathtt{5143267}&\mathtt{4152376}&\mathtt{4376215}&\mathtt{3246517}&\mathtt{5743216}\\
\mathtt{2176534}&\mathtt{2376514}&\mathtt{4652317}&\mathtt{2153467}&\mathtt{3167425}&\mathtt{2174356}\\
\mathtt{2357416}&\mathtt{3165247}&\mathtt{3175246}&\mathtt{3475216}&\mathtt{2475316}&\mathtt{4172365}\\
\mathtt{2564317}&\mathtt{4163275}&\mathtt{2164357}&\mathtt{3564217}&\mathtt{3142576}&
\end{array}
\]
Its 23 pairs have $\det(Z_H^\top Z_H)=155\cdot2^{71}$, and the ratio to the value $2^{78}$ for one pair deleted is $155/128$. After some of the orders are replaced by their reverses, component 1 is in the second position of every order, so the foldover design is obtained by inserting one component into a foldover design with six components and 23 pairs. Read the list row by row and reverse the orders in positions 1, 2, 3, 4, 5, 7, 11, 12, 13, 14, 16 and 19. The resulting 23 orders have $\det(X_D^\top X_D)=355\cdot2^{74}$, and no other choice of one order from each of the 23 pairs gives a larger value. The array of Section~\ref{sec:designs} with one run deleted has $\det(X_D^\top X_D)=2^{82}$, so the ratio is $355/256$.

\noindent\textit{Exchanges of a single run.} For each of the two arrays with five components given above and for the arrays of Examples~3 and~5 of \citet{tsai2022}, every design obtained by deleting one run was examined, together with every exchange of one of its runs for any of the $m!$ orders, repeated runs included. This makes $15\,840$ exchanges for each array with five components, $397\,440$ for six components and $2\,782\,080$ for seven. In integer arithmetic, none of them increases the determinant. Some exchanges leave it unchanged. Putting the deleted run in the place of another run gives the array with that other run deleted. Apart from these exchanges and the exchange of a run for itself, there are 384 exchanges with an unchanged determinant for each of the two arrays with five components, none for the array with six components and 1728 for the array with seven. A sequence of exchanges of this kind can lead to a better design. Delete $13245$ from the first array with five components, and then exchange $43215$ for $12345$, $53214$ for $32145$ and $23415$ for $42315$. The three designs obtained have $|\det X_D|=36\,864$, $36\,864$ and $45\,056$, and each has 11 distinct runs. An exchange algorithm that accepts exchanges with an unchanged determinant can thus leave a design with a run deleted, whereas one that accepts only improvements cannot.

\section{Details of the sensitivity comparison}\label{sec:riskexact}
This section gives the details of Section~6.4 of the paper. Suppose the mean is $\beta_0+z(a)^\top\beta+\phi(a)^\top\gamma$, where the regressors in $\phi(a)\in\mathbb R^r$ have mean zero and are uncorrelated with the PWO factors under the uniform distribution on all orders. Then $\beta$ is the vector of coefficients of the PWO factors in the least squares projection of the mean over the full design, whether or not the additional term is present. Fit the PWO model with an intercept, omitting $\phi$, and put
\[
 P_n=I_n-\mathbf1\mathbf1^\top/n,\qquad
 S_D=Z_D^\top P_nZ_D,\qquad
 \mathcal L_D=S_D^{-1}Z_D^\top P_n.
\]
For independent errors of variance $\sigma^2$, average over random coefficients $\gamma$ with mean zero and $E(\gamma\gamma^\top)=\tau^2I_r/r$, independent of the errors. Write $\Phi_D$ for the matrix with rows $\phi(a)^\top$, $a\in D$. The decomposition into variance and squared bias gives
\begin{equation}\label{eq:studystructuredrisk}
 \begin{split}
 R_D(\tau)
 &:=E\{(\widehat\beta-\beta)^\top\Sigma(\widehat\beta-\beta)\}\\
 &=\sigma^2\operatorname{tr}(\Sigma S_D^{-1})
   +\frac{\tau^2}{r}\operatorname{tr}
       (\Sigma\mathcal L_D\Phi_D\Phi_D^\top\mathcal L_D^\top),
 \end{split}
\end{equation}
because $\widehat\beta-\beta=\mathcal L_D\Phi_D\gamma+\mathcal L_D\varepsilon$ and $\mathcal L_D\mathcal L_D^\top=S_D^{-1}$. Thus $\mathcal V_D=\tr(\Sigma S_D^{-1})$ and $\mathcal B_D=\tr(\Sigma\mathcal L_D\Phi_D\Phi_D^\top\mathcal L_D^\top)/r$ in the notation of the paper. The risk is the mean squared error of the fitted PWO part of the mean, averaged over all orders. It does not include the intercept and does not attempt to predict an unrestricted term $g(a)$ at an order that has not been run. For a foldover design and regressors $\phi$ that are unchanged by reversal, $Z_D^\top\Phi_D=0$ and $Z_D^\top\mathbf1=0$, so the second term is zero. Normality is not needed.

For five components the first family consists of the ten standardised adjacency factors
\[
 \phi^+_{ij}(a)=
 \frac{5\mathbf1_{\{|\pos_a(i)-\pos_a(j)|=1\}}-2}{\sqrt6},
 \qquad 1\leq i<j\leq5.
\]
Each is unchanged by reversal, has mean zero and variance one over the full design, and is uncorrelated with every PWO factor. Here $r=10$. With the stated covariance matrix of $\gamma$, $E_\gamma E_F[(\phi(a)^\top\gamma)^2]=\tau^2$, so $\tau/\sigma$ compares the root mean square of the omitted term with the error standard deviation. The second family consists of the thirty regressors
\[
 \phi^-_{ij,k\ell}(a)=\phi^+_{ij}(a)z_{k\ell}(a),\qquad
 i<j,\ k<\ell,\quad\{i,j\}\cap\{k,\ell\}=\varnothing,
\]
with $r=30$. They change sign under reversal, also have mean zero and variance one, and are uncorrelated with the PWO factors; summation over the 120 orders verifies these properties exactly. Both families and the covariance matrices of their coefficients are unchanged by relabelling of the components, including the changes of sign that this induces. The two families are considered separately.

\begin{figure}[tb]
\centering
\includegraphics[width=\linewidth]{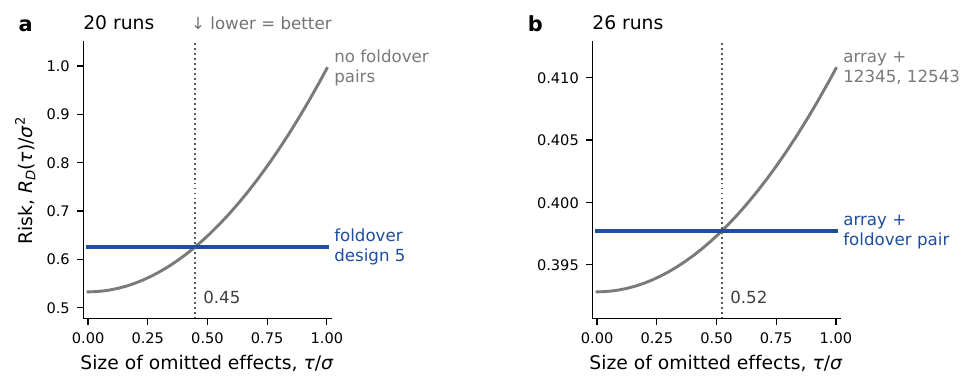}
\caption{Mean squared error $R_D(\tau)/\sigma^2$ of the estimators of the PWO effects when the omitted effects are the adjacency effects, which are unchanged by reversal. Dotted lines mark the values of $\tau/\sigma$ at which the two designs have equal risk. The vertical scales of the two panels differ.}
\label{fig:structuredrisk}
\end{figure}

In $R_D(\tau)/\sigma^2=\mathcal V_D+\mathcal B_D(\tau/\sigma)^2$, the exact coefficients are
\[
\begin{array}{c|ccc}
\text{design}&\mathcal V_D&\mathcal B_D^+&\mathcal B_D^-\\\hline
F_{20}&\frac{5}{8}&0&\frac{103}{72}\\[1ex]
U_{20}&\frac{38909}{73080}&\frac{1481297317}{3204411840}&\frac{663074765}{1922647104}\\[1ex]
F_{26}&\frac{35}{88}&0&\frac{7111}{8712}\\[1ex]
T_{26}&\frac{2885}{7344}&\frac{90475}{5056344}&\frac{91414775}{121352256}
\end{array}
.
\]
Here $F_{20}$ is design 5 of Table~\ref{tab:optimal5}, with half-design
$12345$, $13254$, $14532$, $15423$, $21435$, $24153$, $25314$,
$31524$, $34215$, $43125$; $U_{20}$ is the design without foldover pairs
listed in Section~\ref{sec:enum}; $F_{26}$ is the array generated by~(3.1) of the paper with $12345,54321$ added; and $T_{26}$
is the same array with $12345,12543$ added. The orders $12345$ and $12543$ are not in the array and attain the smallest possible value $|1+\kappa|=1$ for five components, so they maximise the determinant among additions of two runs to this array. For the first family the foldover design has the smaller risk exactly when $\tau/\sigma$ exceeds the value at which the two risks are equal. At 20 and at 26 runs the squares of these values are, respectively,
\[
 (\tau/\sigma)^2=\frac{296675568}{1481297317},\qquad
 (\tau/\sigma)^2=\frac{108783}{398090}.
\]
that is $\tau/\sigma\approx0.447528$ and $0.522745$. Figure~\ref{fig:structuredrisk} shows the two risks at each run size. Adding $12345$ and $12435$, which differ by an interchange of adjacent components, gives a determinant $147/143$ times that of $F_{26}$, and for $T_{26}$ the ratio is $153/143$.

The coefficients above are averages over a family of effects, and for a single omitted effect the comparison can come out differently. Take the regressor $\phi^-_{34,12}$ alone, so that $r=1$, and let $\tau=\sigma$. Among the 25-run designs obtained by adding one order to the array generated by~(3.1) of the paper, the smallest value of $R_D/\sigma^2$ is $23621/21168\approx1.1159$, attained when $12345$ or $54321$ is added. Replacing the run $25341$ of the first of these designs by $31524$ gives a design with 25 distinct runs and $R_D/\sigma^2=44177635/47904048\approx0.9222$. This design has efficiency $0.9728$ under the PWO model, against $0.9935$, and under the average over the thirty regressors of the second family it has $R_D/\sigma^2\approx1.2392$ at $\tau=\sigma$, against $1.2230$. A design of Theorem~1(i) is therefore not the best design for every single omitted effect that changes sign under reversal.

